\documentclass[12pt,reqno]{amsart}

\usepackage{amssymb}
\usepackage{amsthm}
\usepackage{amsmath}
\usepackage{a4wide}
\usepackage{latexsym}
\usepackage[v2,tips]{xy}
\usepackage{mathrsfs}
\usepackage[T1]{fontenc}
\usepackage[latin1]{inputenc}

\usepackage{enumitem}

\newtheorem{theorem}{Theorem}[section]
\newtheorem{lemma}[theorem]{Lemma}
\newtheorem{proposition}[theorem]{Proposition}
\newtheorem{corollary}[theorem]{Corollary}

\theoremstyle{definition}
\newtheorem{definition}[theorem]{Definition}

\newtheorem{example}[theorem]{Example}
\newtheorem{remark}[theorem]{Remark}

\numberwithin{equation}{section}

\newcommand{\N}{\mathbb{N}}

\newcommand{\R}{\mathbb{R}}

\makeatletter
\newcommand\notni{\mathrel{\m@th\mathpalette\canc@l\owns}}
\newcommand\canc@l[2]{{\ooalign{$\hfil#1/\mkern1mu\hfil$\crcr$#1#2$}}}
\makeatother

\newcommand{\smashw}[2][l]{{\text{\makebox[0pt][#1]{$#2$}}}}


\begin{document}
\title[Closed operator ideals
  on Tsirelson and Schreier spaces]{Closed ideals of operators
  on the Tsirelson and Schreier spaces}

\author[K.~Beanland]{Kevin Beanland}
\address[K.~Beanland]{Department of Mathematics, Washington and Lee University, Lexington, VA 24450, USA}
\email{beanlandk@wlu.edu}

\author[T.~Kania]{Tomasz Kania}
\address[T.~Kania]{Mathematical Institute\\Czech Academy of Sciences\\\v Zitn\'a~25 \\115~67 Praha~1\\Czech Republic; and Institute of Mathematics and Computer Science\\ Jagiellonian University\\ {\L}ojasiewicza~6, 30-348~Krak\'{o}w, Poland}
\email{kania@math.cas.cz, tomasz.marcin.kania@gmail.com}

\author[N.~J.~Laustsen]{Niels Jakob Laustsen}
\address[N.~J.~Laustsen]{Department of Mathematics and Statistics, Fylde College, Lancaster University, Lancaster,
LA1 4YF, United Kingdom}
\email{n.laustsen@lancaster.ac.uk}
\begin{abstract}
Let~$\mathscr{B}(X)$ denote the Banach algebra of bounded operators
on~$X$, where~$X$ is either Tsirelson's Banach space or the Schreier
space of or\-der~$n$ for some $n\in\N$. We show that the lattice of
closed ideals of~$\mathscr{B}(X)$ has a very rich structure; in
particular~$\mathscr{B}(X)$ contains at least continuum many maximal
ideals.

Our approach is to study the closed ideals generated by the basis
projections. Indeed, the unit vector basis is an unconditional basis
for each of the above spaces, so there is a~ba\-sis
pro\-jec\-tion~$P_N\in\mathscr{B}(X)$ corresponding to each non-empty
subset~$N$ of~$\N$. A closed ideal of~$\mathscr{B}(X)$ is
\emph{spatial} if it is generated by~$P_N$ for some~$N$. We can now
state our main con\-clu\-sions as follows:
\begin{itemize}
\item the family of spatial ideals lying strictly between the ideal of
  compact operators and~$\mathscr{B}(X)$ is non-empty and has no
  minimal or maximal elements;
\item for each pair $\mathscr{I}\subsetneqq\mathscr{J}$ of spatial
  ideals, there is a family $\{\Gamma_L : L\in
  \Delta\}$, where the index set~$\Delta$ has the cardinality of the
  continuum, such that~$\Gamma_L$ is an uncountable chain of spatial   ideals, $\bigcup\Gamma_L$ is a closed ideal
  that is not spatial, and
\[ \mathscr{I}\subsetneqq\mathscr{L}\subsetneqq\mathscr{J}\qquad\text{and}\qquad \overline{\mathscr{L}+\mathscr{M}} = \mathscr{J} \]
whenever $L,M\in\Delta$ are distinct and $\mathscr{L}\in\Gamma_L$, $\mathscr{M}\in\Gamma_M$.
 \end{itemize}
\end{abstract}
\maketitle
\section{Introduction and statement of main results}
\noindent
Let $X$ be a Banach space with an unconditional basis
$(b_j)_{j\in\N}$.  For a subset~$N$ of~$\N$, we write $P_N$ for the
basis projection corresponding to~$N$; that is, $P_Nx = \sum_{j\in
  N}\langle x,b^*_j\rangle b_j$ for each $x\in X$, where $b^*_j\in X^*$
denotes the $j^{\text{th}}$ coordinate functional.  By a \emph{spatial
  ideal} of the Banach algebra~$\mathscr{B}(X)$ of bounded operators
on~$X$, we understand the closed, two-sided ideal
generated by the basis
projection~$P_N$ for some non-empty subset~$N$ of~$\N$.  A spatial
ideal~$\mathscr{I}$ is \emph{non-trivial} if \[
\mathscr{K}(X)\subsetneqq\mathscr{I}\subsetneqq\mathscr{B}(X), \]
where~$\mathscr{K}(X)$ denotes the ideal of compact operators. A
\emph{chain of spatial ideals} is a non-empty set~$\Gamma$ of spatial
ideals of~$\mathscr{B}(X)$ such that~$\Gamma$ is totally ordered by
inclusion.

We shall study two (classes of) Banach spaces, namely Tsirelson's
space on the one hand and the Schreier spaces of finite order on the
other. We refer to Section~\ref{section:Schreier} for details of the
definition of the latter spaces, originally due to Alspach and
Argyros~\cite{AA}. These Banach spaces have unconditional bases. Using
spatial ideals, we shall show that their lattices of closed operator
ideals have a very rich structure. The following theorem
summarizes our main findings.

\begin{theorem}\label{manychainsofidealsTsirelson}
  Let $X$ denote either Tsirelson's  space~$T$ or the Schreier
  space~$X[\mathcal{S}_n]$ of or\-der~$n$ for some $n\in\N$. 
\begin{enumerate}[label={\normalfont{(\roman*)}}]
\item\label{manychainsofidealsTsirelson1} The family of non-trivial
  spatial ideals of~$\mathscr{B}(X)$ is non-empty and has no minimal or
  maximal elements.
\item\label{manychainsofidealsTsirelson2} Let
  $\mathscr{I}\subsetneqq\mathscr{J}$ be spatial ideals
  of~$\mathscr{B}(X)$. Then there is a family \mbox{$\{\Gamma_L :
    L\in \Delta\}$} such that:
  \begin{itemize}
    \item the index set~$\Delta$ has the cardinality
      of the continuum;
  \item for each $L\in\Delta$, $\Gamma_L$ is an
    uncountable chain of spatial ideals of~$\mathscr{B}(X)$ such that \[
    \mathscr{I}\subsetneqq\mathscr{L}\subsetneqq\mathscr{J}\qquad
    (\mathscr{L}\in\Gamma_L), \]
  and
  $\bigcup\Gamma_L$ is a closed ideal that is not spatial;
  \item $\overline{\mathscr{L}+\mathscr{M}} = \mathscr{J}$ 
 whenever $\mathscr{L}\in\Gamma_L$ and
    $\mathscr{M}\in\Gamma_M$ for  distinct $L,M\in\Delta$.
  \end{itemize}
\item\label{manychainsofidealsTsirelson3} The Banach
  algebra~$\mathscr{B}(X)$ contains at least continuum many maximal
  ideals.
\end{enumerate}
\end{theorem}

We shall also consider the ``small'' ideals of operators on the above
spaces, where we call an ideal ``small'' if it contains no projections
with infinite-dimensional image. The particular ideals that we are
interested in are the compact, strictly singular and inessential
operators; we refer to Definition~\ref{defnSSIness} for the precise definitions of the latter
two operator ideals, which we denote by~$\mathscr{S}$
and~$\mathscr{E}$, respectively. 

\begin{theorem}\label{thmsmallideals}
\begin{enumerate}[label={\normalfont{(\roman*)}}]
\item\label{thmsmallideals1} The ideals of compact, strictly singular
  and inessential operators on Tsirelson's space coincide, and they
  are equal to the intersection of the non-trivial spatial ideals
  of~$\mathscr{B}(T)\colon$
\[ 
\mathscr{K}(T) = \mathscr{S}(T) = \mathscr{E}(T) =
\bigcap\bigl\{\mathscr{I} : \mathscr{I}\ \text{is a non-trivial
  spatial ideal of}\ \mathscr{B}(T)\bigr\}. \]
\item\label{thmsmallideals2} Let $X = X[\mathcal{S}_n]$ be the
  Schreier space of order~$n$ for some $n\in\N$. Then
  \[ \mathscr{K}(X)\subsetneqq\mathscr{S}(X) = \mathscr{E}(X) \]
and
\begin{equation}\label{thmsmallideals2eq} \bigcap\bigl\{\mathscr{I} : \mathscr{I}\ \text{is a non-trivial
  spatial ideal of}\ \mathscr{B}(X)\bigr\}\nsubseteq\mathscr{S}(X). \end{equation}
\end{enumerate}  
\end{theorem}  

The paper is organized as follows: we conclude this introduction with
a brief survey of related results to provide some background for our
work and put it in context. In Sec\-tion~\ref{section:prelim}, we set up
notation and establish some basic general results, as well as a common
frame\-work for the proof of Theorem~\ref{manychainsofidealsTsirelson},
before we complete the proofs for Tsirelson's space in
Sec\-tion~\ref{section:Tsirelson} and for the Schreier spaces in
Sec\-tion~\ref{section:Schreier}. Finally, Section~\ref{sectionOpenQ}
contains some open ques\-tions related to this work.\medskip

The seminal study of operator ideals is due to Calkin~\cite{Calkin},
who considered the situation when the underlying Banach space is a
separable Hilbert space. His most important con\-clu\-sion (at least
from our point of view) is that the ideal of compact operators is the
only proper, non-zero closed ideal in this case. Gohberg, Markus and
Feldman~\cite{GMF} gene\-ra\-lized Calkin's result to the other
classical sequence spaces~$c_0$ and~$\ell_p$ for $p\in\left[1,\infty\right)$,
  while Gramsch~\cite{Gr} and Luft~\cite{luft} independently
  classified all the closed ideals of~$\mathscr{B}(H)$ when~$H$ is a
  non-separable Hilbert space. Their result implies that these ideals
  form a well-ordered chain whose length is determined by the
  dimension of~$H$.
  
Berkson and Porta~\cite[Section~5]{bp} initiated the study of the
closed ideals of~$\mathscr{B}(L_p[0,1])$ for $p\in(1,2)\cup(2,\infty)$,
proving in particular that they are \textsl{not} totally ordered.
Porta~\cite{portaBAMS} then went on to construct a Banach space~$X$
such that there is an injective map from the set~$[\N]^{<\infty}$ of
all finite subsets of the natural numbers into the lattice of closed
ideals of~$\mathscr{B}(X)$, and this map preserves inclusions in both
directions. Porta's Banach space~$X$ is the $\ell_2$-sum of a family
of the form $\{\ell_p : p\in\mathsf{P}\}$ for a countably infinite
subset~$\mathsf{P}$ of~$\left[1,\infty\right)$. By ensuring that
  $1\notin\mathsf{P}$ and that~$\mathsf{P}$ contains the conjugate
  index of each of its elements, Porta arranged that~$X$ is reflexive
  and isometric to its dual space~$X^*$. As far as we know, this was
  the first example of a separable Banach space which has infinitely
  many closed operator ideals.

  Porta~\cite{portaStudia} also initiated the study of the closed operator
  ideals on~$\ell_p\oplus\ell_q$ for distinct
  \mbox{$p,q\in(1,\infty)$}, notably showing that for $p=2$, there are
  exactly two maximal ideals, which correspond to the operators that
  factor through~$\ell_2$ and~$\ell_q$,
  respectively. Subsequently, Volkmann~\cite{vo} extended this result to
  arbitrary finite sums of the
  form $\ell_{p_1}\oplus\ell_{p_2}\oplus\cdots\oplus\ell_{p_m}$
  and $\ell_{p_1}\oplus\ell_{p_2}\oplus\cdots\oplus\ell_{p_{m-1}}\oplus
  c_0$, where $m\in\{2,3,\ldots\}$ and $1\leqslant
  p_1<p_2<\cdots<p_m<\infty$.

Pietsch surveyed these results in his monograph~\cite{pie} and
provided further progress in some cases, observing in particular that there are
infinitely many closed operator ideals on~$L_p[0,1]$ for
$p\in(1,2)\cup(2,\infty)$ and uncountably many on~$C[0,1]$. Moreover,
he formally asked whether there are infinitely many closed operator
ideals on each of the spaces~$L_1[0,1]$, $\ell_p\oplus\ell_q$ and
$\ell_p\oplus c_0$ for $1\leqslant p<q<\infty$. These questions have
only recently been answered, all in the affirmative; we shall give
further details below.

After a relatively quiet period, the study of closed operator ideals
has gained new momentum since the turn of the millenium. Among the
early progress were the first new full classifications of the closed
operator ideals on a Banach space since Gramsch's and Luft's work,
beginning with the space
$\bigl(\bigoplus_{k\in\N}\ell_2^k\bigr)_{c_0}$ and its
dual~$\bigl(\bigoplus_{k\in\N}\ell_2^k\bigr)_{\ell_1}$ (see~\cite{LLR}
and~\cite{lsz}, respectively), and Daws' generalization~\cite{daws} of
Gramsch's and Luft's result to the other non-separable sequence
spaces~$\ell_p(\Gamma)$ and~$c_0(\Gamma)$ for an arbitrary uncountable
index set~$\Gamma$ and $p\in\left[1,\infty\right)$.

Subsequently, as a bi-product of Argyros and Haydon's spectacular
solution~\cite{ah} of the scalar-plus-compact problem, several new
Banach spaces whose closed operator ideals can be classified have
appeared, including~\cite{tarb}, \cite{mpz}
and~\cite{KLindiana}. Another such classification is given
in~\cite[Theorem~5.5]{KK}, namely for the Banach space~$C(K)$,
where~$K$ is Koszmider's Mr\'{o}wka space constructed in~\cite{kosz}
under the Continuum Hypothesis; see~\cite{PKNJL} for the construction of such a space within \textsf{ZFC}. 
We refer to
\cite[Remark~1.5]{KLindiana} for a more detailed survey of the above
results.

An important common feature of the spaces listed in the previous
paragraph is that they are all ``purpose-built'', which is in sharp
contrast to those we described before.  We consider it a very
interesting --- and probably very difficult --- challenge to find new
examples of ``classical'' Banach spaces whose closed operator ideals
can be classified (where ``classical'' can perhaps best be understood
as ``having been known, or at least accessible to Banach''). To
substantiate this claim, we shall outline three further cases, where
apparently ``nice'' Banach spaces have been shown to have very
intricate lattices of closed operator ideals.
Theorem~\ref{manychainsofidealsTsirelson} above in the case of the
original Schreier space~$X[\mathcal{S}_1]$ could arguably also be
included in this list.

We begin with Figiel's reflexive Banach spaces which are not
isomorphic to their cartesian squares~\cite{fi}. These Banach spaces
are manifestly ``nice'', being defined by entirely elementary
means. Indeed, what Figiel showed is that for each strictly decreasing
sequence $(p_k)$ in $(2,\infty)$ and each number $q\in\left(1,\inf p_k\right]$,
there is a sequence $(n_k)$ in~$\N$ such that the Banach space
\[ F =  \Bigl(\bigoplus_{k\in\N}
\ell_{p_k}^{n_k}\Bigr)_{\ell_q} \] satisfies: for each~$m\in\N$,
$F^{m+1}$ does not embed isomorphically into~$F^m$. Now consider the
Boolean algebra $\mathscr{P}(\N)/[\N]^{<\infty}$ which is the quotient
of the power set of~$\N$ modulo the ideal of finite subsets. Loy and
the third-named author \cite[Theorem~4.12]{LR} constructed an
injective map from the ideal lattice
of~$\mathscr{P}(\N)/[\N]^{<\infty}$ into the closed ideal lattice
of~$\mathscr{B}(F)$ such that this map preserves the order  in both
directions. This implies in particular that~$\mathscr{B}(F)$ has
continuum many closed ideals. (The result stated in
\cite[Corollary~4.13]{LR} says ``uncountably many'', but the argument
actually gives continuum many by using the existence of an~al\-most
disjoint family of subsets of~$\N$ having the cardinality of the
continuum.)

More recently, Schlumprecht and Zs\'{a}k~\cite{SchZs} launched the
first successful attack on the above-mentioned questions of Pietsch by
constructing a chain of continuum many closed operator ideals
on~$\ell_p\oplus\ell_q$ for $1<p<q<\infty$. The other cases of this
question have subsequently also been resolved in~\cite{Wa}, \cite{SW}
and~\cite{FSZ}, so we now know that~$\mathscr{B}(\ell_p\oplus\ell_q)$
and $\mathscr{B}(\ell_p\oplus c_0)$ contain uncountable chains of
closed ideals whenever $1\leqslant p<q<\infty$. (In fact, in all cases
except~$\mathscr{B}(\ell_1\oplus c_0)$, the chains have the
cardinality of the continuum.)

Finally, in 2018  Johnson, Pisier and Schechtman~\cite{JPS} answered
the remaining question of Pietsch by constructing a chain of continuum
many closed operator ideals on~$L_1[0,1]$. They also obtained similar
conclusions for~$C[0,1]$ (where previously only uncountably many
closed operator ideals were known) and~$L_\infty[0,1]$ (and therefore also
for~$\ell_\infty$ because~$\ell_\infty$ and~$L_\infty[0,1]$ are isomorphic
as Banach spaces by a theorem of Pe{\l}czy\'{n}ski~\cite{Pel}), using a~variant of their argument and duality, respectively.

\textsl{Note:} after the  initial version of this paper was submitted, Johnson and Schechtman~\cite{JS} have shown that $\mathscr{B}(L_p[0,1])$ contains $2^{\mathfrak{c}}$ closed ideals for each $p\in(1,\infty)\setminus\{2\}$. 

To conclude this survey, let us remark that the Tsirelson and Schreier
spaces are not the first examples of separable Banach spaces having at
least continuum many maximal operator ideals. Indeed,
Mankiewicz~\cite{M} and Dales--Loy--Willis~\cite{DLW} have
independently constructed separable Banach spaces~$X$ such that
$\mathscr{B}(X)$ admits a bounded, surjective algebra
homomorphism~$\varphi$ onto~$\ell_\infty$, and therefore
$\{ \varphi^{-1}(\mathscr{M}) : \mathscr{M}\ \text{is a maximal ideal of}\ \ell_\infty \}$
is a family of cardinality~$2^{\mathfrak{c}}$ of 
maximal ideals of~$\mathscr{B}(X)$. 

\section{Preliminaries, including
the framework of the proof of
Theorem~{\normalfont{\ref{manychainsofidealsTsirelson}}}}\label{section:prelim}
\noindent 
For a set~$N$, $\mathscr{P}(N)$ denotes its power set,
while~$[N]$ and $[N]^{<\infty}$ are the sets of infinite and finite
subsets of~$N$, respectively. We write~$\lvert N\rvert$ for the
cardinality of~$N$; the letter~$\mathfrak{c}$ denotes the cardinality
of the continuum.

For two non-empty subsets $M$ and $N$ of~$\N$, we use the notation
$M<N$ to indicate that~$M$ is finite and $\max M<\min N$.  By an \emph{interval} in a
subset~$N$ of~$\R$, we understand a~set of form~$J\cap N$, where~$J$
is an interval of~$\R$ in the usual sense. (Note that the interval~$J$
may be open, closed or half-open.)

All normed spaces are over the same scalar field~$\mathbb{K}$, which
may be either the real or the complex numbers.  The term ``operator''
means a bounded, linear map between normed spaces.  We
write~$\mathscr{B}(X)$ for the Banach algebra of operators on a Banach 
space~$X$. For $S\in\mathscr{B}(X)$, $\langle
S\rangle$ denotes the (algebraic, two-sided) ideal of~$\mathscr{B}(X)$
generated by~$S$, that is,
\[ \langle S\rangle = \biggl\{\sum_{j=1}^k U_j SV_j : k\in\N,\, U_1,\ldots,U_k,V_1,\ldots,V_k\in\mathscr{B}(X)\biggr\}. \]  
Since $\mathscr{B}(X)$ is a unital Banach algebra, the ideal~$\langle
S\rangle$ is proper if and only if its norm-closure $\overline{\langle
  S\rangle}$ is. The following related result \cite[Lemma~4.9]{LR} is
fundamental to our investigations.

\begin{lemma}\label{LSZobs}
Let~$\mathscr{I}$ be an ideal of a Banach algebra~$\mathscr{A}$, and
let $P\in\mathscr{A}$ be idempotent. Then $P\in\mathscr{I}$ if (and
only if) $P\in\overline{\mathscr{I}}$.
\end{lemma}

\begin{definition}\label{defnSSIness}
Let $X$ and $Y$ be Banach spaces.  An operator $S\colon X\to Y$  is: 
\begin{itemize}
\item \emph{strictly singular} if, for
each $\varepsilon>0$,  each infinite-dimensional
  subspace of~$X$ contains a unit
  vector~$x$ such that $\lVert Sx\rVert<\varepsilon$; in other
  words, the restriction of~$S$ to~$W$ is not an isomorphic embedding
  for any infinite-dimensional
  subspace~$W$ of~$X$; 
\item \emph{inessential} if $I_X + RS$ is a Fredholm operator (that
  is, has finite-dimensional kernel and cofinite-dimensional image)
  for each operator $R\colon Y\to X$, where $I_X$ denotes the identity
  operator on~$X$.  
\end{itemize}  
We write $\mathscr{S}(X,Y)$ and $\mathscr{E}(X,Y)$ for the sets of
strictly singular and inessential operators from~$X$ to~$Y$,
respectively.
\end{definition}
With these definitions, $\mathscr{S}$ and $\mathscr{E}$ define closed
operator ideals in the sense of Pitsch, and
$\mathscr{S}(X,Y)\subseteq\mathscr{E}(X,Y)$ for any Banach spaces~$X$
and~$Y$.  As usual, we write $\mathscr{S}(X)$ and $\mathscr{E}(X)$
instead of $\mathscr{S}(X,X)$ and $\mathscr{E}(X,X)$. A projection
$P\in\mathscr{B}(X)$ is inessential if and only if it has
finite-dimensional image. Although we shall not require this result,
let us mention that $\mathscr{E}(X)$ is equal to the pre-image under
the quotient map of the Jacobson radical of the Calkin algebra
$\mathscr{B}(X)/\mathscr{K}(X)$. This was indeed Kleinecke's original
definition of the inessential operators on a single Banach space~\cite{klein}; the
definition given  above, where the domain and codomain may
differ, is due to Pietsch~\cite{pie}. 

Pfaffenberger~\cite{Pfaf} has shown that $\mathscr{S}(X) =
\mathscr{E}(X)$ whenever the Banach space~$X$ is \emph{sub\-projective}\label{pfaffpageref}
in the sense that each closed, infinite-dimensional subspace
of~$X$ contains a~closed, infinite-dimensional subspace which is
complemented in~$X$. 

Let $X$ and $Y$ be Banach spaces.  A basic sequence $(x_j)_{j\in\N}$
in~$X$ \emph{dominates} a basic sequence $(y_j)_{j\in\N}$ in~$Y$ if
there is a constant $C>0$ such that
\[ \biggl\|\sum_{j=1}^k \alpha_j y_j\biggr\|\leqslant C
\biggl\|\sum_{j=1}^k \alpha_j x_j\biggr\|\qquad (k\in\N,\,
\alpha_1,\ldots,\alpha_k\in\mathbb{K}). \] If we wish to record the
value of the constant~$C$, we say that $(x_j)_{j\in\N}$
\emph{$C$-dominates} $(y_j)_{j\in\N}$.

Let $X$ be a Banach space with an unconditional basis
$(b_j)_{j\in\N}$. It is easy to see that 
the basis projections satisfy the identity
\begin{equation}\label{psisetlemmaNewEq}
 P_{M\cup N} = P_M+P_N - P_{M\cap N}\qquad (M,N\subseteq\N).
\end{equation}
For a  subset~$M$ of~$\N$, we write $X_M$ for the image of the basis projection~$P_M$; that is,
\begin{equation}\label{defnXN} 
X_M = \overline{\operatorname{span}}\{ b_j : j\in M\}. 
\end{equation}
In the notation introduced above, the ideals of the form
$\overline{\langle P_M\rangle}$ for some non-empty subset~$M$ of~$\N$
are precisely the spatial ideals of~$\mathscr{B}(X)$.  The
ideal~$\mathscr{K}(X)$ of compact operators is al\-ways spatial. More
precisely, for $M\subseteq\N$, we have $\overline{\langle
  P_M\rangle}=\mathscr{K}(X)$ if and only if~$M$ is non-empty and
finite.
The following lemma characterizes when one spatial ideal is contained
in another.

\begin{lemma}\label{psisetlemmaNew}
Let $X$ be a Banach space with an unconditional
basis, and let~$M$ and~$N$ be subsets of~$\N$. Then the following four
conditions are equivalent:
\begin{enumerate}[label={\normalfont{(\alph*)}}]
\item\label{psisetlemmaNew2} $P_M\in\overline{\langle P_N\rangle};$
\item\label{psisetlemmaNew1} $\langle P_M\rangle\subseteq\langle
  P_N\rangle;$
\item\label{psisetlemmaNew1.5} $\langle P_N\rangle = \langle
  P_{M\cup N}\rangle;$
\item\label{psisetlemmaNew3}
  $X_M$ is  isomorphic to a complemented subspace of~$X_N^k$
  for some $k\in\N$.
\end{enumerate}
\end{lemma}
\begin{proof}
Lemma~\ref{LSZobs} shows that~\ref{psisetlemmaNew2}
implies~\ref{psisetlemmaNew1}, which in turn
implies~\ref{psisetlemmaNew1.5} by~\eqref{psisetlemmaNewEq}.
 Clearly~\ref{psisetlemmaNew1.5} implies~\ref{psisetlemmaNew2}, and
 finally the equivalence of~\ref{psisetlemmaNew2}
 and~\ref{psisetlemmaNew3} is a special case of \cite[Lemma~4.7]{LLR}
 (or the much earlier \cite[Lemma~1]{portaBAMS} if we know that $X_M\cong X_M\oplus
 X_M$, which will be the case in our applications of this result.)
\end{proof}

\begin{corollary}\label{psisetplusfinitelemma}
 Let $X$ be a Banach space with an unconditional basis, and let~$N$ be
 a~non-empty subset of~$\N$. Then $\langle P_N\rangle = \langle
 P_{N\cup F}\rangle$ for each $F\in[\N]^{<\infty}$.
\end{corollary}

\begin{proof} 
We have $P_F\in\mathscr{F}(X)\subseteq\langle P_N\rangle$ because the
set~$F$ is finite and the ideal~$\mathscr{F}(X)$ of finite-rank operators is the smallest
non-zero ideal of~$\mathscr{B}(X)$. Hence the conclusion
follows from Lemma~\ref{psisetlemmaNew}.
\end{proof}

Combining Lemma~\ref{psisetlemmaNew} with Pe\l{}czy\'{n}ski's Decomposition Method,
we obtain the following conclusion.
\begin{corollary}\label{psisetcorNew}
Let $X$ be a Banach space with an unconditional
basis, and let $M$ and~$N$ be subsets of~$\N$ such that $X_M$ is
isomorphic to~$X_M\oplus X_M$ and~$X_N$ is isomorphic to~$X_N\oplus
X_N$. Then $\langle P_M\rangle = \langle P_N\rangle$ if and only
if~$X_M$ and~$X_N$ are isomorphic.
\end{corollary}

\begin{proposition}\label{propDichotomy} \emph{(Dichotomy for chains
    of spatial ideals)} Let $X$ be a Banach space with an
  unconditional basis, and let~$\Gamma$ be a chain of spatial
  ideals. Then either~$\Gamma$ stabilizes, so that
  $\overline{\bigcup\Gamma}\in\Gamma$, or the ideal
  $\overline{\bigcup\Gamma}$ is not spatial.
\end{proposition}

\begin{proof}
  The two statements are clearly mutually exclusive. Suppose that the
  second statement fails, so that $\overline{\bigcup\Gamma} =
  \overline{\langle P_M\rangle}$ for some non-empty subset~$M$
  of~$\N$. We must show that the first state\-ment is satisfied, that
  is, $\overline{\langle P_M\rangle}\in\Gamma$.  Since a projection
  belongs to the closure of an~ideal if and only if it belongs to the
  ideal itself by Lemma~\ref{LSZobs}, we can find a non-empty
  sub\-set~$N$ of~$\N$ such that $P_M\in\langle P_N\rangle$ and
  $\overline{\langle P_N\rangle}\in\Gamma$. Then
  \[ \overline{\langle P_M\rangle}\subseteq
   \overline{\langle P_N\rangle}\subseteq\overline{\bigcup\Gamma} =
   \overline{\langle P_M\rangle}, \] so we conclude that
   $\overline{\langle P_M\rangle} = \overline{\langle
     P_N\rangle}\in\Gamma$, as required.
\end{proof}  

We shall next state two technical lemmas which will form the core of
the proof of Theorem~\ref{manychainsofidealsTsirelson}. The set-up is
as follows.  Let~$X$ be a Banach space with an unconditional basis,
and suppose that $M\subseteq N$ are non-empty subsets of~$\N$ such
that $P_N\notin\langle P_M\rangle$. We note that~$N$ is infinite
because otherwise $P_N\in\mathscr{F}(X)\subseteq\langle
P_M\rangle$. Further, we see that the set
\[ \Omega_{M,N} = \bigl\{\overline{\langle P_L\rangle} : 
M\subseteq L\subseteq N,\, P_N\notin\overline{\langle
  P_L\rangle}\bigr\} \] of spatial ideals of~$\mathscr{B}(X)$ is
partially ordered by inclusion, and also non-empty with a~smallest
element, namely~$\overline{\langle P_M\rangle}$.  We say that a
chain~$\Gamma$ in~$\Omega_{M,N}$ is \emph{set-induced} if there is
an~in\-creasing sequence $(L_j)_{j\in\N}$ of subsets of~$N$ such that
$M\subseteq L_1$ and $\Gamma = \bigl\{\overline{\langle
  P_{L_j}\rangle} : j\in\N\bigr\}$.
\begin{lemma}\label{newlemma210619}
  Let $X$ be a Banach space with an unconditional basis, and let
  $M\subseteq N$ be non-empty subsets of~$\N$ such that
  $P_N\notin\langle P_M\rangle$.
\begin{enumerate}[label={\normalfont{\Roman*.}}]
\item\label{newlemma210619NewItem0}  A
chain~$\Gamma$ in~$\Omega_{M,N}$ is set-induced if and only if  either~$\Gamma$ stabilizes and has order type~$n$ for some $n\in\N$, or~$\Gamma$ has order type~$\omega$.
\item\label{newlemma210619NewItem1} Suppose that the following two conditions
  are satisfied:
\begin{enumerate}[label={\normalfont{(\ref{newlemma210619NewItem1}\roman*)}}]
\item\label{newlemma210619ii} each set-induced chain in~$\Omega_{M,N}$
  has an upper bound in~$\Omega_{M,N};$
\item\label{newlemma210619i} $\Omega_{M,N}$ has no maximal elements. 
\end{enumerate}
Then each countable chain in $\Omega_{M,N}$ has an upper bound
in~$\Omega_{M,N},$ and there is a~(nec\-es\-sar\-ily uncountable)
chain~$\Gamma$ in $\Omega_{M,N}$ such that:
  \begin{itemize}
  \item $\Gamma$ has no  upper  bound in~$\Omega_{M,N};$
  \item each countable subchain of~$\Gamma$ has an upper bound
    in~$\Gamma;$
  \item the ideal $\bigcup\Gamma$ is closed, and it is not spatial.
  \end{itemize}  
\item\label{newlemma210619NewItem2}  Suppose that there is a map
  $\varphi\colon\mathscr{P}(\N)\to[N]$ which satisfies the following
  three con\-di\-tions for each pair $D,E\in\mathscr{P}(\N)\colon$
  \begin{enumerate}[label={\normalfont{(\ref{newlemma210619NewItem2}\roman*)}}]
  \item\label{newlemma260619i} $M\subseteq \varphi(D);$ 
  \item\label{newlemma260619iii} $\varphi(D)\cup\varphi(E) =
    \varphi(D\cap E);$ 
  \item\label{newlemma260619ii} $P_N\in\langle P_{\varphi(D)}\rangle$
    if and only if $D\in[\N]^{<\infty}.$
  \end{enumerate}
 Then there is a family~$\Delta\subseteq[N]$ of
 car\-di\-nali\-ty~$\mathfrak{c}$ such that
\begin{equation}\label{newlemma260619Eq}
M\subseteq L\quad\text{and}\quad \overline{\langle
  P_M\rangle}\subsetneqq\overline{\langle
  P_L\rangle}\subsetneqq\overline{\langle P_N\rangle} =
\overline{\langle P_{L\cup L'}\rangle}\qquad (L,L'\in\Delta,\,L\ne
L'). \end{equation}
\end{enumerate} 
\end{lemma}

\begin{proof}
  \ref{newlemma210619NewItem0} The forward implication is clear.

  Conversely, let~$\Gamma$ be a chain in~$\Omega_{M,N}$ of order type $n\in\N\cup\{\omega\}$.
  In both cases we can express~$\Gamma$ as $\Gamma =
  \bigl\{\overline{\langle P_{K_j}\rangle} : j\in\N\bigr\}$, where $M\subseteq
  K_j\subseteq N$ and $P_N\notin\overline{\langle
    P_{K_j}\rangle}\subseteq\overline{\langle
    P_{K_{j+1}}\rangle}$
  for each~$j\in\N$. Then, defining
  $L_j = \bigcup_{i=1}^j K_i$,   we obtain
  an increasing sequence~$(L_j)_{j\in\N}$ of subsets of~$N$ such that $M\subseteq L_1$ and
  $\langle P_{K_j}\rangle = \langle P_{L_j}\rangle$ for each $j\in\N$
by Lemma~\ref{psisetlemmaNew}. Con\-sequent\-ly $\Gamma = \{ \overline{\langle P_{L_j}\rangle} : j\in\N\}$ is set-induced.

  \ref{newlemma210619NewItem1}
  Let $\Gamma$ be a countable chain in
  $\Omega_{M,N}$. If~$\Gamma$ is finite, then it has a maximal element and thus an upper bound in~$\Omega_{M,N}$. Otherwise we can choose a subset~$\Upsilon$ of~$\Gamma$ such that~$\Upsilon$ has order type~$\omega$ and every element of~$\Gamma$ is contained in an element of~$\Upsilon$. Then~$\Upsilon$ is set-induced by~\ref{newlemma210619NewItem0} Hence    
  condition~\ref{newlemma210619ii} implies that~$\Upsilon$ has an upper bound
  in~$\Omega_{M,N}$, and that upper bound is clearly also an upper
  bound for~$\Gamma$.

If each chain in~$\Omega_{M,N}$ had an upper bound in~$\Omega_{M,N}$,
then the Kuratowski--Zorn Lemma would imply that~$\Omega_{M,N}$
contains a maximal element, contrary to
condition~\ref{newlemma210619i}.  Therefore~$\Omega_{M,N}$ contains a
chain~$\Gamma$ without any upper bound in~$\Omega_{M,N}$, and this
chain~$\Gamma$ must be
un\-count\-able by the result
proved in the previous paragraph.

To establish the second bullet point, assume towards a contradiction
that~$\Gamma$ contains a~count\-able sub\-chain~$\Upsilon$ which has no
upper bound in~$\Gamma$. As shown above, $\Upsilon$ has an upper bound
$\mathscr{M}\in\Omega_{M,N}$.  The assumption means that, for each
$\mathscr{J}\in\Gamma$, we can find~$\mathscr{L}\in\Upsilon$ such that
$\mathscr{L}\nsubseteq\mathscr{J}$.  Hence
$\mathscr{J}\subseteq\mathscr{L}$ because~$\Gamma$ is a chain, and
therefore also $\mathscr{J}\subseteq\mathscr{M}$. This shows
that~$\mathscr{M}$ is an upper bound for~$\Gamma$, which contradicts
that $\mathscr{M}\in\Omega_{M,N}$.

To see that $\bigcup\Gamma$ is closed, suppose that
$S\in\overline{\bigcup\Gamma}$. We can then recursively construct a
sequence $(S_j)$ of operators and an increasing sequence
$(\mathscr{J}_j)$ of spatial ideals belonging to~$\Gamma$ such that
$S_j\in\mathscr{J}_j$ and $\lVert S- S_j\rVert<1/j$ for each
$j\in\N$. As we showed in the previous paragraph, the countable
subchain $\{\mathscr{J}_j:j\in\N\}$ of~$\Gamma$ has an upper bound
$\mathscr{M}\in\Gamma$. Since ~$\mathscr{M}$ is closed and
contains~$S_j$ for each $j\in\N$, we conclude that
$S\in\mathscr{M}\subseteq\bigcup\Gamma$, as required.  Finally,
since~$\Gamma$ has no upper bound in~$\Omega_{M,N}$, it cannot
stabilize, so Proposition~\ref{propDichotomy} implies that the ideal
$\bigcup\Gamma$ is not spatial.

\ref{newlemma210619NewItem2} We begin by showing that
$\overline{\langle
  P_M\rangle}\subsetneqq\overline{\langle
  P_{\varphi(D)}\rangle}$
whenever $D\subseteq\N$ is
co-infinite. The inclusion follows from~\ref{newlemma260619i}. To see
that it is proper when~$D$ is co-infinite, suppose that
the two ideals are equal, so that $P_{\varphi(D)}\in\langle
P_M\rangle$ by Lemma~\ref{LSZobs}, and set $D^c = \N\setminus
D$. Condition~\ref{newlemma260619i} implies that $P_M\in \langle
P_{\varphi(D^c)}\rangle$, and hence also
$P_{\varphi(D)}\in\langle P_{\varphi(D^c)}\rangle$. Combining this 
with Lemma~\ref{psisetlemmaNew}
and  
\ref{newlemma260619iii}--\ref{newlemma260619ii}, we deduce that
\[ \langle P_{\varphi(D^c)}\rangle = \langle
P_{\varphi(D)\cup \varphi(D^c)}\rangle = \langle
P_{\varphi(\emptyset)}\rangle\ni  P_N. \]
This shows that~$D^c$ is finite by condition~\ref{newlemma260619ii}, and the
conclusion follows.

Now take a family~$\mathscr{D}\subseteq[\N]$ of
cardinality~$\mathfrak{c}$ such that~$\mathscr{D}$ is almost disjoint in the sense that $D\cap
E$ is finite whenever  $D,E\in\mathscr{D}$ are distinct, and set $\Delta =
\varphi(\mathscr{D})$. Then each $D\in\mathscr{D}$ is infinite and
co-infinite, so $\overline{\langle
  P_M\rangle}\subsetneqq\overline{\langle
  P_{\varphi(D)}\rangle}\subsetneqq\overline{\langle P_N\rangle}$ by the
result proved in the previous paragraph and~\ref{newlemma260619ii}. 
Suppose that $D,E\in\mathscr{D}$ are distinct. Then \mbox{$D\cap
E\in[\N]^{<\infty}$},
so $P_N\in\langle P_{\varphi(D\cap E)}\rangle = \langle P_{\varphi(D)\cup
  \varphi(E)}\rangle$ by~\ref{newlemma260619ii} and~\ref{newlemma260619iii}. This establishes the final equality
in~\eqref{newlemma260619Eq}, and it also implies that
$\varphi(D)\ne\varphi(E)$, so~$\varphi$ is injective, and thus
$\lvert\Delta\rvert = \lvert\mathscr{D}\rvert = \mathfrak{c}$.
\end{proof}  

\begin{lemma}\label{keylemma}
 Let $X$ denote either Tsirelson's space~$T$ or the Schreier
 space~$X[\mathcal{S}_n]$ of order~$n$ for some $n\in\N$.
\begin{enumerate}[label={\normalfont{(\roman*)}}]
\item\label{keylemma1} For each  $N\in[\N]$, there is $M\in[N]$ such
  that $P_{N}\notin\langle P_M\rangle$, and consequently
  \[ \overline{\langle P_M\rangle}\subsetneqq\overline{\langle
    P_N\rangle}. \]
\item\label{keylemma2} Suppose that $M\subseteq N$ are infinite
  subsets of~$\N$ such that $P_N\notin\langle P_M\rangle$. Then:
\begin{itemize}
\item   there
  is a map $\varphi\colon\mathscr{P}(\N)\to[N]$ which satisfies
  conditions~\ref{newlemma260619i}--\ref{newlemma260619ii} in
  Lemma~{\normalfont{\ref{newlemma210619},}} and hence there is a 
  family~$\Delta\subseteq[N]$ of car\-di\-nali\-ty~$\mathfrak{c}$ such
  that~\eqref{newlemma260619Eq} is satisfied;
\item each
  set-induced chain in~$\Omega_{M,N}$ has an upper bound
  in~$\Omega_{M,N}$.
\end{itemize}  
\end{enumerate}
\end{lemma}

The proof of Lemma~\ref{keylemma} is non-trivial both for Tsirelson's
space and for the Schreier spaces of finite order; we shall give these
proofs in Sections~\ref{section:Tsirelson} and~\ref{section:Schreier},
respectively.  However, once the lemma is established,
Theorem~\ref{manychainsofidealsTsirelson} follows fairly easily, as we
shall now show.

\begin{proof}[Proof of
    Theorem~{\normalfont{\ref{manychainsofidealsTsirelson},}} assuming
    Lemma~{\normalfont{\ref{keylemma}}}.]  Applying
  Lemma~\ref{keylemma}\ref{keylemma1}--\ref{keylemma2} in the particular case  \mbox{$N=\N$}, we see that $\mathscr{B}(X)$ contains a non-trivial
  spatial ideal and that each proper spatial ideal has at least
  continuum many successors, so no such ideal is maximal.  Another
  application of Lemma~\ref{keylemma}\ref{keylemma1} shows that no
  non-trivial spatial ideal is minimal. This establishes
  Theorem~\ref{manychainsofidealsTsirelson}\ref{manychainsofidealsTsirelson1}.

To verify
Theorem~\ref{manychainsofidealsTsirelson}\ref{manychainsofidealsTsirelson2},
let $\mathscr{I}\subsetneqq\mathscr{J}$ be spatial ideals
of~$\mathscr{B}(X)$, and take non-empty subsets~$K,N$ of~$\N$ such
that $\mathscr{I} = \overline{\langle P_K\rangle}$ and $\mathscr{J} =
\overline{\langle P_N\rangle}$.  By Lemma~\ref{psisetlemmaNew}, we may
replace~$N$ with~$N\cup K$ to ensure that $K\subseteq N$.  Moreover,
we may suppose that~$K$ is infinite. Indeed,
Lemma~\ref{keylemma}\ref{keylemma1} implies that~$N$ contains an
infinite subset~$K'$ such that $\overline{\langle
  P_{K'}\rangle}\subsetneqq \mathscr{J}$, and if~$K$ is finite, then
$\overline{\langle P_{K}\rangle}=\mathscr{K}(X)\subseteq\overline{\langle
  P_{K'}\rangle}$, so we may replace~$K$ with~$K'$.

This enables us to apply Lemma~\ref{keylemma}\ref{keylemma2} with $M = K$ to
obtain a family $\Delta\subseteq[N]$ of
car\-di\-nali\-ty~$\mathfrak{c}$ such that
\begin{equation}\label{keylemma2EqK}
K\subseteq L\quad\text{and}\quad
\mathscr{I}\subsetneqq\overline{\langle
  P_L\rangle}\subsetneqq\mathscr{J} = \overline{\langle P_{L\cup
    L'}\rangle}\qquad (L,L'\in\Delta,\,L\ne L'), \end{equation} and
each set-induced chain in~$\Omega_{K,N}$ has an upper bound
in~$\Omega_{K,N}$.  Take $L\in\Delta$. Then \mbox{$P_{N}\notin\langle
  P_L\rangle$}, and Lemma~\ref{keylemma}\ref{keylemma2} (this time
applied with $M=L$) shows that the pair $L\subseteq N$ satisfies
conditions~\ref{newlemma210619ii}--\ref{newlemma210619i} in
Lemma~\ref{newlemma210619}.  Consequently~$\Omega_{L,N}$ contains an
uncountable chain~$\Gamma_L$ such that $\bigcup\Gamma_L$ is a closed
ideal that is not spatial. It follows that the first two bullet points
in
Theorem~\ref{manychainsofidealsTsirelson}\ref{manychainsofidealsTsirelson2}
are satisfied. To verify the third, suppose that
$\mathscr{L}\in\Gamma_L$ and $\mathscr{M}\in\Gamma_{L'}$, where
$L,L'\in\Delta$ are distinct.  Then $P_L\in\mathscr{L}$ and
$P_{L'}\in\mathscr{M}$, and therefore~\eqref{psisetlemmaNewEq} shows
that $P_{L\cup L'}\in\mathscr{L}+\mathscr{M}$. Hence the conclusion
follows from the fact that $\overline{\langle P_{L\cup
    L'}\rangle}=\mathscr{J}$
by~\eqref{keylemma2EqK}.

\ref{manychainsofidealsTsirelson3}. Applying
clause~\ref{manychainsofidealsTsirelson2} 
in the case where $\mathscr{J} =\mathscr{B}(X)$ and~$\mathscr{I}$ is
any proper spatial ideal, we deduce that there is a family~$\Xi$ of
proper spatial ideals of~$\mathscr{B}(X)$ such that~$\Xi$ has
cardinality~$\mathfrak{c}$ and
\begin{equation}\label{manychainsofidealsTsirelson3Eq}
\overline{\mathscr{L}+\mathscr{M}}=\mathscr{B}(X)\qquad
(\mathscr{L},\mathscr{M}\in\Xi,\, \mathscr{L}\neq\mathscr{M}).
\end{equation}
Each of the ideals in~$\Xi$ is contained in a maximal ideal
of~$\mathscr{B}(X)$, and~\eqref{manychainsofidealsTsirelson3Eq}
implies that these maximal ideals are all distinct,
so~$\mathscr{B}(X)$ contains at least continuum many maximal ideals.
\end{proof} 

\section{Tsirelson's space}\label{section:Tsirelson}
\noindent Following Figiel and Johnson~\cite{FJ}, we use the term
\emph{Tsirelson's space} for the dual of the reflexive Banach space
that Tsirelson~\cite{T} originally constructed with the property that
it does not contain any of the classical sequence spaces~$c_0$
and~$\ell_p$ for $1\leqslant p<\infty$, and we denote it by~$T$. This
convention makes no difference from the point of view of ideal
lattices because~$T$ is reflexive, so the adjoint map $S\mapsto
S^*,\,\mathscr{B}(T)\to\mathscr{B}(T^*),$ is an isometric, linear
bijection which is anti-multiplicative in the sense that
$(RS)^*=S^*R^*$ for $R,S\in\mathscr{B}(T)$, and therefore it induces a
lattice isomorphism between the closed ideal lattices
of~$\mathscr{B}(T)$ and~$\mathscr{B}(T^*)$.

We refer to Casazza and Shura's monograph~\cite{CS} for details about
Tsirelson's space, including its formal definition. In line with their notation, we
write~$(t_j)_{j\in\N}$ for the unit vector basis, which is a
normalized, $1$\nobreakdash-un\-condi\-tional basis for~$T$. Recall
from~\eqref{defnXN} that~$T_M$ denotes the closed linear span of $\{
t_j : j\in M\}$ in~$T$ for a subset~$M$ of~$\N$. Using this no\-ta\-tion,
we have the following fundamental result \cite[Corollary VII.b.3]{CS}. 
\begin{theorem}\label{thmCSVIIb3old}
Let $M,N\in[\N]$. Then $T_M$ is isomorphic to~$T_N$ if and only if
$(t_j)_{j\in M}$ is equivalent to~$(t_j)_{j\in N}$.
\end{theorem} 

The usefulness of this result relies on being able to determine when
two subsequences of the basis $(t_j)$ are equivalent. Fortunately,
Casazza, Johnson and Tzafriri~\cite{CJT} have identified an index which
does exactly that. To define it, we require some notation. First, for
$M = \{m_1<m_2<\cdots\}\in[\N]$ and $J\in[\N]^{<\infty}$, let
$\sigma(M,J)$ denote the norm of the formal identity operator from the linear span of~$(t_{m_j})_{j\in J}$ to~$\ell_1(J)$, that is,
\[ \sigma(M,J) = \sup\biggl\{ \sum_{j\in J} \alpha_j : \alpha_j\in[0,1],\, \Bigl\|\sum_{j\in J}\alpha_j t_{m_j}\Bigr\|_T\leqslant1\biggr\}, \]
with the convention that $\sigma(M,\emptyset) =0$.  Second, suppose
that \mbox{$N = \{n_1<n_2<\cdots\}\in[\N]$}, and set $m_0 =
n_0=0$. Then Casazza, Johnson and Tzafriri
 have shown that~$(t_j)_{j\in M}$ is
equivalent to~$(t_j)_{j\in N}$ if and only if
\begin{equation}\label{CJTcrit} \sup\Bigl\{\sigma\bigl(M, M\cap\left(n_{j-1},n_{j}\right]\bigr), \sigma\bigl(N, N\cap\left(m_{j-1},m_{j}\right]\bigr): j\in\N\Bigr\}<\infty \end{equation}
 (see \cite[the remark following Theorem~10]{CJT}).  This result
simplifies considerably in the special case where $M\subseteq N$,
which will suffice for our purposes. We incorporate it in the
following omnibus characterization of equality of spatial ideals of~$\mathscr{B}(T)$,
which will be our key tool in the proof of Lemma~\ref{keylemma} for Tsirelson's space.

\begin{corollary}\label{containmentofspatialideals}
The following conditions are equivalent for each pair $M\subseteq N$ of infinite subsets of~$\N\colon$
\begin{enumerate}[label={\normalfont{(\alph*)}}]
\item\label{containmentofspatialideals1}
 $P_N\in \overline{\langle P_M\rangle};$
\item\label{containmentofspatialideals2}
$\langle P_M\rangle = \langle P_N\rangle;$
\item\label{containmentofspatialideals3}
$T_N$ is  isomorphic to a complemented subspace of~$T_M;$
\item\label{containmentofspatialideals4}
$T_N$ is isomorphic to~$T_M;$
\item\label{containmentofspatialideals5} $(t_j)_{j\in M}$ is equivalent to  $(t_j)_{j\in N};$
\item\label{containmentofspatialideals6} there is a constant $C\geqslant 1$ such that $\sigma(N,J)\leqslant C$
  for each interval~$J$ in~$N$ with $J\cap M=\emptyset$.
\end{enumerate}
\end{corollary}
\noindent 
As it will be used repeatedly in the remainder of this section, let us
spell out that the conditions on the set~$J$ in
clause~\ref{containmentofspatialideals6} above mean that $J =
N\cap[a,b]$ for some numbers $a, b$ satisfying $m_{j-1}<a\leqslant
b<m_{j}$ for some $j\in\N$, where $m_0=0$ and $M= \{m_1<m_2<\cdots\}$
is the increasing enumeration of~$M$, as above.

\begin{proof}
The assumption that $M\subseteq N$ means that $P_M=P_MP_N\in\langle
P_N\rangle$, and hence
conditions~\ref{containmentofspatialideals1}--\ref{containmentofspatialideals4}
are equivalent by Lemma~\ref{psisetlemmaNew} and
Corollary~\ref{psisetcorNew} because $T_L\cong T_L\oplus T_L$ for each
$L\in[\N]$, as explained in \cite[the paragraph following Proposition~I.12]{CS}. Theorem~\ref{thmCSVIIb3old} shows that
conditions~\ref{containmentofspatialideals4}
and~\ref{containmentofspatialideals5} are equivalent.

Finally, conditions~\ref{containmentofspatialideals5}
and~\ref{containmentofspatialideals6} are equivalent by the result of
Casazza, Johnson and Tzafriri stated above.  Indeed, suppose
that~\ref{containmentofspatialideals5} is satisfied, and
let~$C\in\left[1,\infty\right)$ be the supremum given by~\eqref{CJTcrit}. Since each interval~$J$ in~$N$ with $J\cap M=\emptyset$ is contained in $N\cap (m_{j-1},m_j)$ for some $j\in\N$, we have
  $\sigma(N,J)\leqslant \sigma\bigl(N,N\cap
  \left(m_{j-1},m_j\right]\bigr)\leqslant C$.

Conversely, suppose that~$C$ is a constant such that~\ref{containmentofspatialideals6} is satisfied, and take $j\in\N$. Since $M\subseteq N$, the set $M\cap\left(n_{j-1},n_j\right]$ is either empty or a singleton, so $\sigma\bigl(M,M\cap\left(n_{j-1},n_j\right]\bigr)\leqslant 1$. Moreover, the subadditivity of the operator norm implies that
\[ \sigma\bigl(N,N\cap \left(m_{j-1},m_j\right]\bigr)\leqslant \sigma\bigl(N,N\cap (m_{j-1},m_j)\bigr)+1\leqslant C+1, \] and  hence the supremum in~\eqref{CJTcrit} is at
most~$C+1$.
\end{proof}

\begin{proof}[Proof of Lemma~{\normalfont{\ref{keylemma}}} for $X=T$.]
  \ref{keylemma1}.  Let $N= \{n_1<n_2<\cdots\}\in[\N]$. For each
  $k\in\N$, we can find $m>k$ such that
  \mbox{$\sigma\bigl(N,N\cap(k,m)\bigr)> k$} because otherwise the
  basic sequence $(t_{n_j})_{j\in N\cap(k,\infty)}$ would
  $k$-dominate, and hence be equivalent to, the unit vector basis
  of~$\ell_1$ for some $k\in\N$. Using this observation, we can
  recursively construct a strictly increasing sequence $(m_j)$ in~$N$
  such that $\sigma\bigl(N,N\cap(m_{j-1},m_{j})\bigr)>m_{j-1}$ for each
  $j\in\N$, where $m_0 = 0$. Now
  Corollary~\ref{containmentofspatialideals} shows that the subset $M
  = \{m_1<m_2<\cdots\}$ of~$N$ has the desired property.

  \ref{keylemma2}. Let $M\subseteq N$ be infinite subsets of~$\N$ such
  that $P_N\notin\langle P_M\rangle$.  By
  Corollary~\ref{containmentofspatialideals}, we can recursively
  choose intervals $J_1<J_2<\cdots$ in~$N$ such that $J_k\cap
  M=\emptyset$ and $\sigma(N,J_k)\geqslant k$ for each $k\in\N$.
  We shall show that the map $\varphi\colon\mathscr{P}(\N)\to[N]$
defined by 
\begin{equation*}
\varphi(D) =
N\setminus\bigcup_{j\in D} J_j
\end{equation*}
 satisfies conditions~\ref{newlemma260619i}--\ref{newlemma260619ii} in
 Lemma~\ref{newlemma210619}. The first two of these conditions are
 immediate. To verify the third, suppose that $D\subseteq\N$ is
 finite. Then the set $\bigcup_{j\in D} J_j$ is also finite, and
 therefore $P_N\in\langle P_{\varphi(D)}\rangle$ by
 Corollary~\ref{psisetplusfinitelemma}. Conversely, suppose that
 $D\subseteq\N$ is infinite.  For each $k\in D$, $J_k$ is an interval
 in~$N$ such that $J_k\cap \varphi(D)=\emptyset$ and
 $\sigma(N,J_k)\geqslant k$, so as~$D$ is unbounded,
 Corollary~\ref{containmentofspatialideals} implies that
 $P_N\notin\langle P_{\varphi(D)}\rangle$.  This
 establishes~\ref{newlemma260619ii} and hence completes the proof of
 the first bullet point.

To verify the second, let~$\Gamma$ be a set-induced chain
in~$\Omega_{M,N}$, say $\Gamma = \bigl\{\overline{\langle
  P_{L_j}\rangle} : j\in\N\bigr\}$, where $M\subseteq L_j\subseteq
L_{j+1}\subseteq N$ and $P_N\notin\overline{\langle P_{L_j}\rangle}$
for each $j\in\N$.  By Corollary~\ref{containmentofspatialideals}, we
may recursively construct intervals $J_1<J_2<\cdots$ in~$N$ such that
$\sigma(N,J_j)\geqslant j$ and $J_j\cap L_j=\emptyset$ for each
$j\in\N$.  Set $K = N\setminus\bigcup_{j\in\N} J_j$. Then $M\subseteq
K\subseteq N$, and $P_N\notin\overline{\langle P_K\rangle}$ by
Corollary~\ref{containmentofspatialideals} because~$J_j$ is an
interval in~$N$ with $K\cap J_j=\emptyset$ and $\sigma(N,J_j)\geqslant
j$ for each $j\in\N$. Hence $\overline{\langle P_K\rangle}$ belongs
to~$\Omega_{M,N}$. Moreover, for each $k\in\N$ and $j\geqslant k$, we
have $J_j\cap L_k\subseteq J_j\cap L_j=\emptyset$, so $L_k\subseteq
K\cup\bigcup_{j<k} J_j$. Since $\bigcup_{j<k} J_j$ is finite,
Corollary~\ref{psisetplusfinitelemma} implies that $P_{L_k}\in\langle
P_K\rangle$, and therefore~$\overline{\langle P_K\rangle}$ is an upper
bound for~$\Gamma$.
\end{proof} 

\enlargethispage*{10pt}
We require the following result, which is \cite[Proposition~II.7]{CS},
 to prove
Theorem~\ref{thmsmallideals}\ref{thmsmallideals1}. 
\begin{theorem}\label{CSpropII7}
Every closed, infinite-dimensional subspace of~$T$ contains a closed
subspace which is complemented in~$T$ and isomorphic to~$T_N$ for some
$N\in[\N]$.
\end{theorem}

\begin{proof}[Proof of
    Theorem~{\normalfont{\ref{thmsmallideals}\ref{thmsmallideals1}}}]
  A standard perturbation argument shows that $\mathscr{K}(T) =
  \mathscr{S}(T)$, as remarked in~\cite[p.~1173]{ABK}, for instance.
  As observed in \cite[Proposition~2.4(5)]{GMASB},
  Theorem~\ref{CSpropII7} implies that~$T$ is subprojective, and
  therefore $\mathscr{S}(T) = \mathscr{E}(T)$ by the
  result~\cite{Pfaf} of Pfaffenberger mentioned on
  page~\pageref{pfaffpageref}.

The inclusion  
\begin{equation}\label{Eq:14072019} \mathscr{K}(T)\subseteq\bigcap\bigl\{\mathscr{I} : \mathscr{I}\ \text{is a non-trivial spatial ideal of}\ \mathscr{B}(T)\bigr\} \end{equation}
is clear. Conversely, suppose that
$S\in\mathscr{B}(T)\setminus\mathscr{K}(T)$. Then, as explained in the
first paragraph, $S$ is not strictly singular. Take a closed,
infinite-dimensional subspace~$W$ of~$T$ such that the restriction
of~$S$ to~$W$ is an isomorphism onto its image~$S(W)$. Then~$S(W)$ is
a closed, infinite-dimensional subspace of~$T$, so
Theorem~\ref{CSpropII7} implies that~$S(W)$ contains a closed
subspace~$Z$ which is complemented in~$T$ and isomorphic to~$T_N$ for
some $N\in[\N]$. Let $Q\colon T\to Z$ be a projection, and let
$U\colon Z\to T_N$ be an isomorphism.  We observe that the restriction
of~$S$ to the subspace $Y = S^{-1}(Z)\cap W$ is an isomorphism
onto~$Z$; denote it by $\widetilde{S}$. Then we have a commutative
diagram
 \begin{equation*}
    \spreaddiagramrows{2ex}\spreaddiagramcolumns{3ex}%
    \xymatrix{T\ar[rrrrrr]^{\displaystyle{P_N}}\ar@{-->}[dd]
      \ar[drr]^{\displaystyle{\widetilde{S}^{-1}U^{-1}P_N}} & & & & &
      & T\\ & &
      Y\ar[rr]^{\displaystyle{\widetilde{S}}}_{{\cong}}\ar@{_{(}->}[dll]
      & & Z\ar[r]^{\displaystyle{U}}_{{\cong}} &
      T_N\ar@{^{(}->}[ur]\\ T\ar[rrrrrr]^{\displaystyle{S}} & & & & &
      &
      T\smashw{,}\ar@{-->}[uu]\ar[ull]_{\displaystyle{Q}}} \end{equation*}
 where the two unlabelled solid arrows are the set-theoretic
 inclusions. This diagram shows that~$P_N$ factors through~$S$, and
 therefore $P_N\in\langle S\rangle$. By
 Lemma~\ref{keylemma}\ref{keylemma1}, we can find $M\in[N]$ such that
 $P_N\notin\overline{\langle P_M\rangle}$. Consequently
 $S\notin\overline{\langle P_M\rangle}$, so~$S$ does not belong to the
 right-hand side of~\eqref{Eq:14072019}.
\end{proof}

\section{The Schreier spaces of finite order}\label{section:Schreier}

\noindent
The aim of this section is to establish Lemma~\ref{keylemma} and
Theorem~\ref{thmsmallideals}\ref{thmsmallideals2} for the Schreier
space~$X[\mathcal{S}_n]$ of order~$n\in\mathbb{N}$ associated with the
Schreier family~$\mathcal{S}_n$, originally defined by Alspach and
Argyros~\cite{AA}. Their precise definition is as follows.

\begin{definition}  Let \[ \mathcal{S}_0 = \bigl\{\{k\} : k\in\mathbb{N}\bigr\}\cup\{\emptyset\}, \]
and for $n\in\mathbb{N}_0$, recursively define
\[ \mathcal{S}_{n+1} = \biggl\{\bigcup_{i=1}^k E_i : k\in\mathbb{N},\, E_1,\ldots,E_k\in\mathcal{S}_n\setminus\{\emptyset\},\,k\leqslant \min E_1,\, E_1<E_2<\cdots<E_k\biggr\}\cup\{\emptyset\}.\]
For $n\in\mathbb{N}_0$, the \emph{Schreier space of order}~$n$ is the
completion of~$c_{00}$ with respect to the norm
\[ \|x\| = \sup\biggl\{\sum_{j\in E}|\alpha_j| : E\in\mathcal{S}_n\setminus\{\emptyset\}\biggr\}\qquad \bigl(x=(\alpha_j)_{j\in\N}\in c_{00}\bigr). \]
We denote this Banach space by $X[\mathcal{S}_n]$.
\end{definition}

Of course, the Schreier space of order~$0$ is simply~$c_0$.  For a
fixed order $n\in\N_0$, we write $(e_j)_{j\in\N}$ for the unit
vector basis for~$X[\mathcal{S}_n]$. Alspach and Argyros~\cite{AA}
have shown that this basis is $1$-un\-condi\-tional and shrink\-ing.

\textsl{Note:} the original definition of Alspach and Argyros of the
Schreier family~$\mathcal{S}_n$ has an ex\-ten\-sion to the case where~$n$
is a countably infinite ordinal. We have stated it for finite~$n$ only because we have been unable to extend our
results beyond that case.

Each subset of a set in~$\mathcal{S}_n$ is clearly also
in~$\mathcal{S}_n$. Another elementary and often useful property of
the Schreier family~$\mathcal{S}_n$ is that it is \emph{spreading} in
the following sense. Let $J = \{ j_1<j_2<\cdots<j_m\}$ and \mbox{$K
  =\{k_1<k_2<\cdots<k_m\}$} be finite subsets of~$\N$, and
suppose that~$K$ is a \emph{spread} of~$J$; that is,
$j_i\leqslant k_i$ for each $i\in\{1,\ldots,m\}$. Then
$J\in\mathcal{S}_n$ implies that $K\in\mathcal{S}_n$.

We shall require several results and definitions from the paper~\cite{GaspLeung} of Gasparis and Leung that we shall now review. However, the starting point is a result of Gasparis alone \cite[Corollary~3.2]{Gasp}, which states that, for each $n\in\N_0$ and $M\in [\mathbb{N}]$, there is a unique sequence
$\bigl(F_j^n(M)\bigr)_{j\in\N}$ of finite subsets of~$M$ such that: 
\begin{enumerate}\label{GaspResult}
\item $M= \bigcup_{j\in\N} F_j^n(M);$
\item $F_j^n(M)$ is a maximal $\mathcal{S}_n$-set for each $j\in\N$ in
  the sense that $E=F_j^n(M)$ is the only set $E\in \mathcal{S}_n$ such
  that $F_j^n(M)\subseteq E;$
\item the sets $F_j^n(M)$ are successive in the sense that $F_j^n(M)<F_{j+1}^n(M)$ for each $j\in\N$.
\end{enumerate}
Gasparis and Leung~\cite[Defini\-tion~3.1]{GaspLeung} used this result
to define the following numerical 
index for a set
$J\in[\mathbb{N}]^{<\infty}$:
\[ \tau_n(J)=\max\{k\in\N : J \cap F_k^n(J_{\text{tail}})\ne\emptyset\}, \]
where we have introduced the notation
$J_{\text{tail}} = J\cup \{ j\in\N : j>\max J\}$ in an attempt to make
the expression a little easier to comprehend.
Roughly speaking, $\tau_n(J)$ counts how many successive maximal
$\mathcal{S}_n$-sets (almost) fit inside~$J$.  The following
remark collects three easy observations concerning this index.
\begin{remark}\label{tauncharRemark} Let $n\in\N_0$ and $J\in [\mathbb{N}]^{<\infty}$. Then:
\begin{enumerate}[label={\normalfont{(\roman*)}}]
\item\label{tauncharRemark1} $\tau_n(J)\geqslant 1$ if and only if
  $J\ne\emptyset.$
\item\label{tauncharRemark3} For $k\in\N$, $\tau_n(J)=k$ if and only
  if there are successive $\mathcal{S}_n$-sets $E_1<E_2<\cdots< E_k$
  such that $J= \bigcup_{i=1}^k E_i$ and $E_1,\ldots,E_{k-1}$ are
  maximal $\mathcal{S}_n$-sets. (Note that the final
  $\mathcal{S}_n$-set~$E_k$ need not be maximal.)
\item\label{tauncharRemark2}
Suppose that~$J$ is non-empty. Then
\begin{equation*}
 \tau_n(J) = \min\Bigl\{ k\in\N : J\subseteq\bigcup_{i=1}^k
 E_i,\ \text{where}\ E_1,\ldots,E_k\in\mathcal{S}_n\ \text{and}\ E_1<E_2<\cdots<E_k\Bigr\}.
\end{equation*}
\end{enumerate}
\end{remark}

As Gasparis and Leung observed \cite[Lemma~3.2(2)]{GaspLeung}, the latter
formula implies that~$\tau_n$ is sub\-addi\-tive in the sense that
\begin{equation}\label{gaspleung322} \tau_n(J\cup
  K)\leqslant \tau_n(J)+\tau_n(K)
\end{equation}
whenever $J,K\in[\N]^{<\infty}$ are successive.  

Gasparis and Leung used the index~$\tau_n$ to define another
index~$d_n$, which can be viewed as a (not necessarily symmetric) way of
measuring the distance from one infinite subset of~$\N$ to another in
terms of the Schreier family~$\mathcal{S}_n$. To help state the
definition of~$d_n$ clearly and compactly, we introduce the following
piece of notation, which was not used by Gasparis and Leung: let~$J$
and~$K$ be (finite or infinite) subsets of~$\N$ such that $\sup
J\leqslant |K|$ (so that in particular~$K$ is infinite whenever~$J$ is
infinite), and enumerate~$K$ in increasing order: $K =
\{k_1<k_2<\cdots\}$.  Then we set
\begin{equation}\label{LJdefn} K(J) = \{ k_j : j\in J\}. \end{equation}

\begin{definition}[Gasparis and Leung {\cite[Defini\-tion~3.3]{GaspLeung}}]\label{GaspLeungDefn3.3}
Let $n\in\N_0$ and $M,N\in [\mathbb{N}]$. The $n^{\text{th}}$ \emph{Gasparis--Leung index} of~$M$ with respect to~$N$ is given by
\[ d_n(M,N)=\sup \bigl\{\tau_n\bigl(M(J)\bigr): J\in[\N]^{<\infty},\,
N(J)\in\mathcal{S}_n\bigr\}. \]
\end{definition}

We note that $d_n(M,N) = 1$ whenever $M\subseteq N$ because
in this case~$M(J)$ is a spread of~$N(J)$ for each
$J\in[\N]^{<\infty}$. 

The significance of the Gasparis--Leung index~$d_n$ is due to the
following result (see \cite[Corollary~1.2(1) and Lemma~3.4, including
  its proof]{GaspLeung}), which together with the immediate
consequence that we record in Corollary~\ref{CorGLthm} will be a key
tool for us. Recall that the notation~$X_N$ was introduced
in~\eqref{defnXN}.
\begin{theorem}[Gasparis and Leung]\label{Thmbasisdomination}
Let $X=X[\mathcal{S}_n]$ for some $n\in\N$, and let $M,N\in [\mathbb{N}]$. 
\begin{enumerate}[label={\normalfont{(\roman*)}}]
\item\label{Thmbasisdomination1} Suppose that~$X_M$ is isomorphic to a
  subspace of~$X_N$. Then the basic sequence $(e_j)_{j\in M}$
  dominates~$(e_j)_{j\in N}.$
\item\label{Thmbasisdomination2} The basic sequence $(e_j)_{j\in M}$
  dominates~$(e_j)_{j\in N}$ if and only if $d_n(M,N)$ is
  finite.\\ Moreover, when $d_n(M,N)$ is finite, $(e_j)_{j\in M}$
  $d_n(M,N)$-dominates $(e_j)_{j\in N}$.
\end{enumerate}
\end{theorem}

\begin{corollary}\label{CorGLthm}
  Let $X=X[\mathcal{S}_n]$ for some $n\in\N$, and let $M,N\in
  [\mathbb{N}]$. Then the following three conditions are equivalent:  
  \begin{enumerate}[label={\normalfont{(\alph*)}}]
  \item\label{CorGLthm1} the subspaces~$X_M$ and~$X_N$ are isomorphic;
  \item\label{CorGLthm2}  $X_M$ is isomorphic to a subspace
    of~$X_N,$ and~$X_N$ is isomorphic to a subspace of~$X_M;$
 \item\label{CorGLthm3} $d_n(M,N)$ and $d_n(N,M)$ are both finite.
\end{enumerate}
\end{corollary}  

Our first application of these results is
to establish the following proposition.
\begin{proposition}\label{XMcartesian} Let $X = X[\mathcal{S}_n]$ for some $n\in\N$. Then~$X_M$ is isomorphic to~$X_M\oplus X_M$ for each $M\in[\N]$.
\end{proposition}

For clarity of presentation, we have split the proof into a series of lemmas.
\begin{lemma}\label{leftshiftoperatorbounded} Let $n\in\N_0$, and let  $\sigma\colon\N\to\N$ be a strictly increasing map. Then the left shift
\[ L_\sigma\colon\ \sum_j \alpha_j e_j\mapsto \sum_j \alpha_{\sigma(j)} e_j,\quad c_{00}\to c_{00}, \]
extends to an operator of norm one on~$X[\mathcal{S}_n]$.
\end{lemma}
\begin{proof}
For $x=(\alpha_j)_{j\in\N}\in c_{00}$, choose $E\in\mathcal{S}_n\setminus\{\emptyset\}$
such that $\|L_\sigma x\| = \sum_{j\in E}|\alpha_{\sigma(j)}|$.
The set~$\sigma(E)$ is a spread of~$E$ because~$\sigma$ is strictly
increasing. Hence $\sigma(E)\in\mathcal{S}_n$, and therefore
\[ \| x\|\geqslant\sum_{k\in\sigma(E)}|\alpha_k| = \sum_{j\in E}|\alpha_{\sigma(j)}| =  \|L_\sigma x\|. \qedhere \]
\end{proof}

\begin{lemma}\label{lemmaKB20Aug} Let $n\in\N$, and let $E <F$ be successive
  maximal $\mathcal{S}_n$-sets. Then:
\begin{enumerate}[label={\normalfont{(\roman*)}}]
\item\label{lemmaKB20Aug1}
  $|E\cup F|\geqslant 3\min E$.
\item\label{lemmaKB20Aug2} Suppose that $\min E\geqslant 2$. Then $2(E
  \cup F)\notin\mathcal{S}_n$.
\end{enumerate}  
\end{lemma}

\begin{proof} \ref{lemmaKB20Aug1}.  The maximality of~$E$ and~$F$ means that  $|E|\geqslant\min E$ and \[ |F|\geqslant\min F\geqslant \max
  E+1\geqslant\min E + |E|\geqslant 2\min E, \] and hence we have
  $|E\cup F|=|E|+|F|\geqslant 3\min E$.

  \ref{lemmaKB20Aug2}. We proceed by induction on~$n$. 

  The result follows easily from~\ref{lemmaKB20Aug1} for $n=1$ because
  in this case we have
  \[   |2(E\cup F)| = |E\cup F|\geqslant 3\min E>2\min E=\min2(E\cup
  F), \] so that $2(E \cup F)\notin\mathcal{S}_1$.

  Now assume inductively that the result holds for some $n\in\N$. To
  prove it for $n+1$, let $E <F$ be maximal $\mathcal{S}_{n+1}$-sets
  with $k:=\min E\geqslant 2$, and set $\ell=\min F$. Then we can find
  maximal $\mathcal{S}_n$-sets $E_1<\cdots<
  E_k<E_{k+1}<\cdots<E_{k+\ell}$ such that $E=\bigcup_{i=1}^k E_i$ and
  $F= \bigcup_{i=1}^\ell E_{k+i}$. Applying~\ref{lemmaKB20Aug1} to the
  maximal $\mathcal{S}_n$-sets $E_1<E_2$, both of which are contained
  in~$E$, we obtain $|E|\geqslant|E_1\cup E_2|\geqslant 3\min E_1=
  3k$, and therefore
  \[ 4k\leqslant \min E + |E|\leqslant \max E + 1\leqslant \ell.\] 
  In particular $4k\leqslant k+\ell$, so that 
  \[ \bigcup_{i=1}^{2k} 2(E_{2i-1}\cup E_{2i})\subseteq 2(E\cup F), \]
  where $2(E_1\cup E_2),2(E_3\cup E_4),\ldots,2(E_{4k-1}\cup E_{4k})$
  are successive sets which do not belong to~$\mathcal{S}_n$ by the
  induction hypothesis. Hence, if $2(E\cup F)$ is written as the union
  of~$m$ successive $\mathcal{S}_n$-sets for some $m\in\N$, then
  we must have $m>2k = \min 2(E\cup F)$. This shows that $2(E\cup
  F)\notin\mathcal{S}_{n+1}$, and hence the induction continues.
\end{proof}

\begin{corollary}\label{cor25aug}
For each $n\in\N$, $d_n(\N,2\N)\leqslant 3$, and hence 
the right shift given by
\begin{equation}\label{cor25augEq}
 R\colon\ \sum_j \alpha_j e_j\mapsto \sum_j \alpha_j e_{2j},\quad
 c_{00}\to c_{00}, \end{equation} extends to an operator of norm at
most~$3$ on~$X[\mathcal{S}_n]$.
\end{corollary}

\begin{proof} Using~\eqref{LJdefn}
  and Definition~\ref{GaspLeungDefn3.3}, we see that $d_n(\N,2\N)\leqslant 3$
  if and only if $\tau_n(J)\leqslant 3$ for each $J\in[\N]^{<\infty}$
  with $2J\in\mathcal{S}_n$.  By contraposition, the latter statement
  is equivalent to the statement that $2J\notin\mathcal{S}_n$ for each
  $J\in[\N]^{<\infty}$ with $\tau_n(J)\geqslant 4$. To verify this,
  suppose that $J\in[\N]^{<\infty}$ with $\tau_n(J)\geqslant
  4$. Since~$\{1\}$ is a maximal $\mathcal{S}_n$-set, we have
  $\tau_n(J\setminus\{1\})\geqslant 3$, so by
  Remark~\ref{tauncharRemark}\ref{tauncharRemark3}, there are maximal
  $\mathcal{S}_n$-sets $E <F$ such that $E\cup F\subseteq
  J\setminus\{1\}$. Lemma~\ref{lemmaKB20Aug}\ref{lemmaKB20Aug2} then
  implies that $2(E \cup F)$ does not belong to~$\mathcal{S}_n$, and
  the same is therefore true for its superset~$2J$, as desired.

  To establish that the right shift~$R$ given by~\eqref{cor25augEq} is
  bounded by~$3$, we simply combine the inequality
  $d_n(\N,2\N)\leqslant 3$ with
  Theorem~\ref{Thmbasisdomination}\ref{Thmbasisdomination2} to deduce
  that
\[ \biggl\|\sum_{j=1}^k\alpha_j e_{2j}\biggr\|\leqslant 3\biggl\|\sum_{j=1}^k\alpha_j e_{j}\biggr\|\qquad (k\in\N,\,\alpha_1,\ldots,\alpha_k\in\mathbb{K}). \qedhere \]
\end{proof} 

The following example shows that we cannot in general lower the upper bound~$3$ on the quantity  $d_n(\N,2\N)$ in the above proof. It also shows that it is possible to have $\| R\|>2$. 
\begin{example} Let $n=2$, and consider the set $J= \{1,2,\ldots,8\}$. We see that $\tau_2(J) = 3$ because $J = \{1\}\cup\{2,3,4,5,6,7\}\cup\{8\}$, where the first two sets on the right-hand side are maximal $\mathcal{S}_2$-sets. However, $2J$ belongs to~$\mathcal{S}_2$ because it is the union of the two $\mathcal{S}_1$-sets $\{2,4\}$ and $\{6,8,10,12,14,16\}$. Hence $d_2(\N,2\N)\geqslant 3$.  

Set $x = e_1 + \frac{1}{6}(e_2+e_3+\cdots + e_8)\in X[\mathcal{S}_2]$. Then the above reasoning shows that~$\|x\| =1$ (attained at $\{1\}$ and at any $6$-element subset of $\{2,3,\ldots,8\}$), but $\| Rx\| =  1 + \frac{7}{6} >2$ because its support belongs to~$\mathcal{S}_2$. 
\end{example}

\begin{lemma}\label{dnMM'bound} Let  $n\in\N$ and 
  $M\in[\N]$, and set $M' = (2M-1)\cup (2M)\in[\N]$.  Then
  \begin{equation}\label{dnMM'boundEq1} d_n(M,M')\leqslant 3\qquad\text{and}\qquad d_n(M',M)\leqslant
  2, \end{equation} and hence $X_M\cong X_{M'}$, where $X =
  X[\mathcal{S}_n]$.
\end{lemma}
\begin{proof}
Write $M=\{m_1<m_2<\cdots\}$ and $M'=\{m_1'<m_2'<\cdots\}$, where
$m_{2j-1}' = 2m_j-1$ and $m_{2j}' = 2m_j$ for each $j\in\N$.

To prove the first inequality in~\eqref{dnMM'boundEq1}, suppose that
$J\in[\N]^{<\infty}$ with $M'(J)\in\mathcal{S}_n$. The definitions
above imply that $m_j'\leqslant 2m_j$ for each $j\in\N$, so that
$2M(J)$ is a spread of $M'(J)$ and therefore
$2M(J)\in\mathcal{S}_n$. Hence $\tau_n(M(J))\leqslant 3$ by
Corollary~\ref{cor25aug}, and the conclusion follows from
Definition~\ref{GaspLeungDefn3.3}.

Interchanging the roles of~$M$ and~$M'$, we see that the second
inequality in~\eqref{dnMM'boundEq1} amounts to showing that
$\tau_n(M'(J))\leqslant 2$ for each non-empty set $J\in[\N]^{<\infty}$
with $M(J)\in\mathcal{S}_n$. We shall establish this estimate by
induction on~$n\in\N_0$.

Note that we start the induction at $n=0$ for convenience. Indeed,
the estimate is clear in this case because the non-empty
$\mathcal{S}_0$-sets are precisely the singletons, so in fact
$M(J)\in\mathcal{S}_0$ implies that $M'(J)\in\mathcal{S}_0$.

Now assume inductively that we have established the estimate for some
$n\in\N_0$, and let $J\in[\N]^{<\infty}$ be a non-empty set with
$M(J)\in\mathcal{S}_{n+1}$.  Set $j = \min J$. Then, by the definition
of~$\mathcal{S}_{n+1}$, we can find $h\in\N$ and $J_1<J_2<\cdots<J_h$
such that $h\leqslant m_j$, $J = \bigcup_{i=1}^h J_i$ and
$M(J_i)\in\mathcal{S}_n$ for each $i\in\{1,\ldots,h\}$.

Take $k\in\N$ and $L_1<L_2<\cdots<L_k$ such that $J = \bigcup_{i=1}^k
L_i$, $M'(L_i)$ is a maximal $\mathcal{S}_n$-set for each
$i\in\{1,\ldots,k-1\}$ and $M'(L_k)\in\mathcal{S}_n$. The induction
hypothesis implies that $k\leqslant 2h$ because, by
Remark~\ref{tauncharRemark}\ref{tauncharRemark3}, each of the
sets~$J_i$ can be split into at most two successive pieces
\mbox{$J_i'<J_i''$} with $M'(J_i'), M'(J_i'')\in\mathcal{S}_n$.  If
$k\leqslant m_j' = \min M'(J)$, then $M'(J)\in\mathcal{S}_{n+1}$, and
the conclusion follows, so we may suppose that $k>m_j'$. We observe
that $m_j'\geqslant j$ because $m_j'$ is the $j^{\text{th}}$ element
of a strictly increasing sequence of natural numbers, and hence $E :=
\bigcup_{i=1}^j M'(L_i)\in\mathcal{S}_{n+1}$. Since the sets
$L_1,\ldots,L_{j+1}$ are successive and non-empty, we see that $\min
L_{j+1}\geqslant\min L_1+j = 2j$, so that $\min M'(L_{j+1})\geqslant
m_{2j}' = 2m_j\geqslant 2h\geqslant k$, which implies that $F :=
\bigcup_{i=j+1}^k M'(L_i)\in\mathcal{S}_{n+1}$. This shows that
$\tau_{n+1}(M'(J))\leqslant 2$ because $M'(J) = E\cup F$, and hence
the induction continues.

The final clause is immediate from Corollary~\ref{CorGLthm}. 
\end{proof}
\enlargethispage*{.1pt}
\begin{proof}[Proof of Proposition~{\normalfont{\ref{XMcartesian}}}]
  Let $M\in[\N]$, and set $M' = (2M-1)\cup (2M)\in[\N]$.  Then
  $X_M\cong X_{M'}$ by Lemma~\ref{dnMM'bound}, and $X_{M'}\cong
  X_{2M-1}\oplus X_{2M}$ because $M'$ is the disjoint union of the
  sets~$2M-1$ and~$2M$. Hence the result will follow provided that we
  can show that $X_{2M-1}$ and~$X_{2M}$ are both isomorphic to~$X_M$.
  
  The map $\sigma\colon j\mapsto 2j,\,\N\to\N,$ is strictly
  increasing, and the corresponding left shift $L_\sigma$ is a left
  inverse of the right shift~$R$, using the notation of
  Lemma~\ref{leftshiftoperatorbounded} and
  Corollary~\ref{cor25aug}. Hence the restriction of~$R$ to~$X_M$ is
  an isomorphism onto its image, which is~$X_{2M}$, with the inverse
  being the appropriate restriction of~$L_\sigma$.

  A similar argument using the strictly increasing map $\sigma\colon
  j\mapsto 2j-1,\,\N\to\N,$ and the right shift given by
  $\sum_j\alpha_je_j\mapsto \sum_j\alpha_je_{2j-1}$ (which is bounded
  because it equals the composition $L_\rho\circ R$, where $\rho\colon
  j\mapsto j+1,\,\N\to\N,$ and $R$ is given by~\eqref{cor25augEq} as above)
  shows that $X_M\cong X_{2M-1}$.
\end{proof}

Using these results, we obtain the following characterization of when
two spatial ideals of~$\mathscr{B}(X[\mathcal{S}_n])$ are
equal. It is the counterpart of
Corollary~\ref{containmentofspatialideals} and will play a similar
role in our proof of Lemma~\ref{keylemma} for the Schreier spaces of
finite order.

\begin{proposition}\label{ThmPsiSetStrong}
  Let $X=X[\mathcal{S}_n]$ for some $n\in\N$, and suppose
  that~$M,N\in[\N]$ satisfy $P_M\in\overline{\langle
    P_N\rangle}$. Then the following conditions are equivalent:
\begin{enumerate}[label={\normalfont{(\alph*)}}]
\item\label{ThmPsiSetStrong0} $P_N\in\overline{\langle P_M\rangle};$
\item\label{ThmPsiSetStrong1} $\langle P_M\rangle = \langle P_N\rangle;$
\item\label{ThmPsiSetStrong2} $X_M$ is isomorphic to~$X_N;$
\item\label{ThmPsiSetStrong3} $X_N$ is  isomorphic to a subspace of~$X_M;$
\item\label{ThmPsiSetStrong4} the $n^{\text{th}}$ Gasparis--Leung index $d_n(N,M)$ is finite;
\item\label{ThmPsiSetStrong5} there is a constant $k\in\N$ such that
  $\tau_n(N(J))\leqslant k$ for each
  set~$J\in[\N\cap(k,\infty)]^{<\infty}$ with $M(J)\in\mathcal{S}_n$.
\end{enumerate}
\end{proposition}

\begin{proof} Conditions~\ref{ThmPsiSetStrong0} and~\ref{ThmPsiSetStrong1} are equivalent by Lemma~\ref{LSZobs}, while
  Proposition~\ref{XMcartesian} implies that
  Corollary~\ref{psisetcorNew} applies, and therefore
  conditions~\ref{ThmPsiSetStrong1} and~\ref{ThmPsiSetStrong2} are
  also equivalent.  

Condition~\ref{ThmPsiSetStrong2} trivially implies~\ref{ThmPsiSetStrong3}, which in turn implies~\ref{ThmPsiSetStrong4} by Theorem~\ref{Thmbasisdomination}. 
Combining the assumption that $P_M\in\overline{\langle P_N\rangle}$
with Lemma~\ref{psisetlemmaNew} and Proposition~\ref{XMcartesian}, we
deduce that~$X_N$ contains a complemented subspace which is isomorphic
to~$X_M$, and therefore $d_n(M,N)$ is finite by
Theorem~\ref{Thmbasisdomination}. Thus the implication
\ref{ThmPsiSetStrong4}$\Rightarrow$\ref{ThmPsiSetStrong2} follows from
Corollary~\ref{CorGLthm}.

 Definition~\ref{GaspLeungDefn3.3} shows that~\ref{ThmPsiSetStrong4}
 implies~\ref{ThmPsiSetStrong5}.  Conversely, suppose that
 \mbox{$k\in\N$} is a constant such that~\ref{ThmPsiSetStrong5} is
 satisfied. Then, for each set $J\in[\N]^{<\infty}$ with
 $M(J)\in\mathcal{S}_n$, the subadditivity of~$\tau_n$ stated
 in~\eqref{gaspleung322} implies that
  \[ \tau_n(N(J))\leqslant \tau_n\bigl(N(J\cap[1,k])\bigr) + \tau_n\bigl(N(J\cap(k,\infty))\bigr)\leqslant 
  \tau_n\bigl(N(\{1,\ldots,k\})\bigr) + k, \] and therefore
  $d_n(N,M)\leqslant \tau_n\bigl(N(\{1,\ldots,k\})\bigr) + k<\infty$.
\end{proof} 

\begin{lemma}\label{LemmaSpread}
  Let $M,N\in[\mathbb{N}]$, and suppose that $J\in[\N]^{<\infty}$ is a
  non-empty set such that $N\cap \left[1,\max M(J)\right)\subseteq
    M$. Then~$N(J)$ is a spread of~$M(J)$. 
\end{lemma}

\begin{proof}
  Write $M=\{m_1<m_2<\cdots\}$ and $N=\{n_1<n_2<\cdots\}$, and let
  $j\in J$. We must show that $m_j\leqslant n_j$. This is clear if
  $n_j\geqslant\max M(J)$ because $m_j\in M(J)$. Otherwise $n_j\in
  N\cap [1,\max M(J))$, which is contained in~$M$ by the assumption,
    so that $n_j = m_k$ for some $k\in\N$. We have $k\geqslant j$
    because $n_1,\ldots,n_{j-1}\in\{ m_1,\ldots,m_{k-1}\}$, and
    therefore $m_j\leqslant m_k = n_j$, as required.
\end{proof}

\begin{proof}[Proof of Lemma~{\normalfont{\ref{keylemma}}} for {$X=X[\mathcal{S}_n],\,n\in\N$}.]
\ref{keylemma1}.  Let $N\in[\N]$, and let
$\bigl(F_j^n(N)\bigr)_{j\in\N}$ be the unique sequence of successive,
maximal $\mathcal{S}_n$-sets partitioning~$N$ described on
page~\pageref{GaspResult}. The fact that the sets
$F_1^n(N),F_2^n(N),\ldots$ are successive and partition~$N$ means that
we can partition~$\N$ into successive intervals $J_1<J_2<\cdots$ such
that $F_j^n(N)=N(J_j)$ for each $j\in\N$, where $N(J_j)$ (unlike
$F_j^n(N)$) is defined using the notation~\eqref{LJdefn}.  Set $k_1 =
1$ and recursively define $k_{j+1} = k_j+j$ for $j\in\N$. (In other
words, $k_{j+1} = j(j+1)/2+1$, but this formula is not helpful for our
purposes). Then, setting $K_j = \bigcup_{i=k_j}^{k_{j+1}-1}
J_i\in[\N]^{<\infty}$ for each $j\in\N$, we obtain a partition of~$\N$
into successive intervals such that
\[ N(K_j) = \bigcup_{i=k_j}^{k_{j+1}-1} F_i^n(N). \]
Thus~$N(K_j)$ is the union of $k_{j+1}-k_j = j$ successive, maximal
$\mathcal{S}_n$-sets, so $\tau_n(N(K_j)) = j$.

Since $n\geqslant 1$, $N$ contains arbitrarily long
$\mathcal{S}_n$-sets. This fact enables us to recursively choose
successive intervals $L_1<L_2<\cdots$ in~$\N$ such that $\lvert
L_j\rvert = \lvert K_j\rvert$ and $N(L_j)\in\mathcal{S}_n$ for each
$j\in\N$. Indeed, once the intervals $L_1<\cdots <L_j$ have been
chosen for some $j\in\N$, we can take $\ell>\max L_j$ so large that
the interval $L_{j+1}=\left[\ell,\ell+\lvert K_j\rvert\right)\cap\N$ satisfies
  $N(L_{j+1})\in\mathcal{S}_n$, and hence the recursion continues.

Set $L = \bigcup_{j\in\N}L_j\in[\N]$, and observe that $L(K_j)=L_j$
for each $j\in\N$ because $(K_j)$ is a partition of~$\N$ into
successive intervals with $\lvert K_j\rvert = \lvert
L_j\rvert$. Consequently $M := N(L)\in[N]$ satisfies
\[ M(K_j) = (N(L))(K_j) = N(L(K_j)) = N(L_j)\in\mathcal{S}_n, \] 
so $d_n(N,M)\geqslant \tau_n(N(K_j))=j$.  Since this is true for every
$j\in\N$, Proposition~\ref{ThmPsiSetStrong} implies that
$P_N\notin\langle P_M\rangle$.

\ref{keylemma2}. Let $M\subseteq N$ be infinite subsets of~$\N$ with
$P_N\notin\langle P_M\rangle$. By
Corollary~\ref{psisetplusfinitelemma}, we may suppose that $\min M =
\min N$ by adding the element $\min N$ to the set~$M$ if necessary. In
order to define a map~$\varphi\colon\mathscr{P}(\N)\to[N]$ which satisfies con\-di\-tions~\ref{newlemma260619i}--\ref{newlemma260619ii} in Lemma~\ref{newlemma210619}, we shall construct a sequence $(J_i)_{i\in\N}$ of
finite, successive intervals of~$\N$ such that:
\begin{enumerate}[label={\normalfont{(\roman*)}}]
  \item\label{defnIntervalPartition2} $\bigcup_{i\in\N} J_i = \N$; 
  \item\label{defnIntervalPartition3} for each $j\in\N$, $J_j$
    contains a subset~$K_j$ such that $\tau_n(N(K_j))\geqslant j$ and the
    set 
\begin{equation}\label{defnMj}
L_j = M\cup\bigcup_{i<j} N(J_i)\in[N]
\end{equation}
satisfies $L_j(K_j)\in\mathcal{S}_n$ and $L_j(K_j)\subseteq N(J_j)$.
\end{enumerate}
The construction is by recursion, where
condition~\ref{defnIntervalPartition2} is replaced with the
appropriate finite analogue, that is,
\begin{enumerate}[label={\normalfont{(\roman*$'$)}}]
\item\label{defnIntervalPartition2'} $\bigcup_{i=1}^j J_i = \N\cap
  [1,\max J_j]$ for each $j\in\N$.
\end{enumerate}

To begin the recursion, we define
$J_1=K_1=\{1\}$. Then~\ref{defnIntervalPartition2'} is obvious,
and~\ref{defnIntervalPartition3} follows almost as easily because
$L_1=M$ by definition and $M(K_1) = \{\min M\}= N(J_1)$; being a~singleton, this set belongs to~$\mathcal{S}_n$.

Now assume recursively that, for some $j\in\N$, we have chosen finite,
successive inter\-vals \mbox{$J_1<\cdots<J_j$} such that
conditions~\ref{defnIntervalPartition2'}
and~\ref{defnIntervalPartition3} are
satisfied. Following~\eqref{defnMj}, we define 
\mbox{$L_{j+1} =
M\cup\bigcup_{i\leqslant j} N(J_i)$}.  Then
$L_{j+1}\setminus M$ is finite, so
Corollary~\ref{psisetplusfinitelemma} implies that \[ \langle
  P_{L_{j+1}}\rangle = \langle P_{M}\rangle\notni P_N. \]
Hence, by Prop\-o\-si\-tion~\ref{ThmPsiSetStrong}, we can find a set
$K_{j+1}\in[\N\cap(\max J_j,\infty)]^{<\infty}$ such that
\mbox{$L_{j+1}(K_{j+1})\in\mathcal{S}_n$} and $\tau_n(N(K_{j+1}))>
\max J_j$.  We see that $\max J_j\geqslant j$ because the sets
$J_1,\ldots,J_j$ are non-empty and successive, and consequently
$\tau_n(N(K_{j+1}))\geqslant j+1$. Since $L_{j+1}\subseteq N$, we can
choose $k\in\N$ such that $\max L_{j+1}(K_{j+1})$ is the
$k^{\text{th}}$ element of~$N$. Note that \mbox{$k\geqslant\max K_{j+1}$}, so
that $J_{j+1} = \N\cap\left(\max J_j,k\right]$ is a finite successor interval
  of~$J_j$ con\-tain\-ing~$K_{j+1}$ and such
  that~\ref{defnIntervalPartition2'} is satisfied for $j+1$. Moreover,
  \mbox{$\max L_{j+1}(K_{j+1})\in N(J_{j+1})$} by the choice of~$k$, while
  \[ L_{j+1}(K_{j+1})> L_{j+1}\Bigl(\bigcup_{i\leqslant j} J_i\Bigr)=
  N\Bigl(\bigcup_{i\leqslant j} J_i\Bigr), \] so
  $L_{j+1}(K_{j+1})\subseteq N(J_{j+1})$ because~$N(J_{j+1})$ is the
  immediate successor interval 
  of $N\bigl(\bigcup_{i\leqslant
    j}\! J_i\bigr)$ in~$N$. This shows that~\ref{defnIntervalPartition3} is
  also satisfied for $j+1$, and hence the recursion continues.

For $D\subseteq\N$, set $D^c = \N\setminus D$, and define
$\varphi\colon\mathscr{P}(\N)\to[N]$ by 
\begin{equation}\label{SchreierContinuumSuccessorsPsiEq1}
  \varphi(D)=M\cup\bigcup_{i\in D^c}N(J_i).
\end{equation}
This map clearly satisfies
conditions~\ref{newlemma260619i}--\ref{newlemma260619iii} in
Lemma~\ref{newlemma210619}. To help us establish
con\-di\-tion~\ref{newlemma260619ii}, we shall show that
\begin{equation}\label{SchreierContinuumSuccessorsPsiEq2}
  \varphi(D)\cap\left[1,\max L_j(K_j)\right)\subseteq L_j\qquad (D\in\mathscr{P}(\N),\,j\in D),
\end{equation}
where $L_j$ is given by~\eqref{defnMj}. Indeed, suppose that $k\in
\varphi(D)\cap\left[1,\max L_j(K_j)\right)$ for some $j\in D\subseteq\N$. Since
  $M\subseteq L_j$, it suffices to consider the case where $k\in
  N(J_i)$ for some $i\in D^c$. We must have $i\leqslant j$ because \[
  k<\max L_j(K_j)\leqslant\max N(J_j). \] Moreover, $i\ne j$ because
  $i\in D^c$ and $j\in D$. Hence $i<j$, so $L_j\supseteq N(J_i)\ni k$,
  as desired.

Now let $D\subseteq\N$. Suppose first that $P_N\in\langle
P_{\varphi(D)}\rangle$, and let $j\in D$.  Combining
Lemma~\ref{LemmaSpread}
with~\eqref{SchreierContinuumSuccessorsPsiEq2}, we see that
$(\varphi(D))(K_j)$ is a spread of~$L_j(K_j)$, so that
$(\varphi(D))(K_j)\in\mathcal{S}_n$, and therefore $d_n(N,
\varphi(D))\geqslant \tau_n(N(K_j))\geqslant j$. Thus the set~$D$ is
bounded above by $d_n(N, \varphi(D))$, which is finite by
Proposition~\ref{ThmPsiSetStrong}. Hence~$D$ is finite.  Conversely,
suppose that~$D$ is finite. Then $N\setminus \varphi(D)\subseteq
\bigcup_{i\in D} N(J_i)$ is also finite, and therefore $P_N\in\langle
P_{\varphi(D)}\rangle$ by Corollary~\ref{psisetplusfinitelemma}.  This
completes the proof of~\ref{newlemma260619ii} and hence of the first
bullet point in Lemma~\ref{keylemma}\ref{keylemma2}.

To establish the second, let $\Gamma = \bigl\{\overline{\langle
  P_{L_j}\rangle} : j\in\N\bigr\}$ be a set-induced chain
in~$\Omega_{M,N}$, where $M\subseteq L_j\subseteq L_{j+1}\subseteq N$
and $P_N\notin\overline{\langle P_{L_j}\rangle}$ for each $j\in\N$.
We shall recursively choose non-empty, finite subsets $J_1,J_2,\ldots$
of~$\N$ such that
\begin{equation}\label{SchreierupperboundsforcountablechainsEq1}
\tau_n(N(J_j))\geqslant j,\qquad
L_j(J_j)\in\mathcal{S}_n\qquad\text{and}\qquad
L_{j+1}(J_{j+1})>L_j(J_j)\qquad (j\in\N). \end{equation} We begin this
recursion by taking $J_1=\{1\}$. Now assume that non-empty, finite
subsets $J_1,\ldots,J_j$ of~$\N$ have been chosen for some $j\in\N$,
and take $k\in\N$ such that $\max L_j(J_j)$ is the $k^{\text{th}}$
element of~$L_{j+1}$. Since $P_N\notin\overline{\langle
  P_{L_{j+1}}\rangle}$, Proposition~\ref{ThmPsiSetStrong} enables us
to choose a set $J_{j+1}\in[\N\cap(k,\infty)]^{<\infty}$ such that
$\tau_n(N(J_{j+1}))>j$ and $L_{j+1}(J_{j+1})\in\mathcal{S}_n$. Then
the first two statements
in~\eqref{SchreierupperboundsforcountablechainsEq1} are satisfied for
$j+1$, while the last part follows from the fact that
\[ L_j(J_j)\subseteq
L_{j+1}(\N\cap[1,k])< L_{j+1}(J_{j+1}). \] Hence the recursion
continues.

Set $\ell_0 = 0$ and $\ell_j = \max L_j(J_j)$ for $j\in\N$, and define
$L = \bigcup_{j\in\N}L_j\cap \left(\ell_{j-1},\ell_j\right]$. We check
that $\overline{\langle P_L\rangle}\in\Omega_{M,N}$:
\begin{itemize}
\item each $m\in M$ belongs to $\left(\ell_{j-1},\ell_j\right]$ for some
  $j\in\N$, and also $m\in L_j$ because $M\subseteq L_j$, so $m\in L;$
\item $L\subseteq N$ because $L_j\subseteq N$ for each $j\in\N;$
\item for each $j\in\N$, we have
  \[ L\cap [1,\ell_j] = \bigcup_{i=1}^j L_i\cap\left(\ell_{i-1},\ell_i\right]\subseteq
  L_j, \] so Lemma~\ref{LemmaSpread} implies that $L(J_j)$ is a spread
  of~$L_j(J_j)$, and therefore $L(J_j)\in\mathcal{S}_n$; hence
  $d_n(N,L)\geqslant \tau_n(N(J_j))\geqslant j$, so as $j\in\N$ was
  arbitrary, we conclude that $P_N\notin\overline{\langle P_L\rangle}$
  by Proposition~\ref{ThmPsiSetStrong}.
\end{itemize}
We observe that $L_j\cap(\ell_{j-1},\infty)\subseteq L$ for each
$j\in\N$ because the sequence $(L_j)$ is increasing.  This implies
that $P_{L_j}\in\langle P_L\rangle$ by
Corollary~\ref{psisetplusfinitelemma}, and thus $\overline{\langle
  P_L\rangle}$ is an upper bound for~$\Gamma$, as desired.
\end{proof}

\begin{lemma}\label{formalidXtoc0}  For each $n\in\N$, the formal identity operator from~$X[\mathcal{S}_n]$ to~$c_0$ is a
  strictly singular, non-compact operator of norm one.
\end{lemma}

\begin{proof} Set $X=X[\mathcal{S}_n]$, and let $S\colon
  \operatorname{span}\{ e_j : j\in\N\}\to c_0$ be the formal identity
  operator, that is, the linear map determined by $S e_j = d_j$ for
  each $j\in\N$, where $(d_j)$ denotes the unit vector basis
  for~$c_0$. Writing $(e_j^*)$ for the coordinate functionals in~$X^*$
  corresponding to the basis~$(e_j)$ for~$X$, we observe that
  \begin{equation}\label{formalidXtoc0Eq1} \lVert S x\rVert_\infty = \sup\bigl\{ \lvert\langle
    x,e_j^*\rangle\rvert : j\in\N\bigr\} \end{equation} for each
  $x\in\operatorname{span}\{ e_j : j\in\N\}$, and therefore~$S$ is
  bounded with norm~$1$, so it extends uniquely to an operator defined
  on all of~$X$, also denoted by~$S$ and still of norm~$1$; for later
  reference, we note that~\eqref{formalidXtoc0Eq1} remains valid for
  each $x\in X$. This operator~$S$ cannot be compact because~$(e_j)$
  is a bounded sequence in~$X$ such that no subsequence of~$(Se_j)$ converges in~$c_0$.
  
  Assume towards a contradiction that~$S$ is not strictly singular. Then~$X$ contains a closed,
  infinite-dimensional subspace~$W$ such that there exists an
  $\varepsilon>0$ for which \mbox{$\lVert S w\rVert_\infty\geqslant
  \varepsilon \lVert w\rVert_X$} for each $w\in W$.  Choose $m\in\N\cap
  \bigl[2(1+\varepsilon)/\varepsilon^2,\infty\bigr)$, and set $k_1 =
    m$. We can then recursively choose numbers $k_2,\ldots,k_{m+1}$
    with $k_{j+1}>k_j$ and unit vectors $w_j\in
    W\cap\overline{\operatorname{span}}\{ e_i : i\geqslant k_j\}$ such
    that $\bigl|\langle w_j, e_i^*\rangle\bigr|\leqslant\varepsilon/m$
    for each $i\geqslant k_{j+1}$ and $j\in\{1,\ldots,m\}$.

  Set $w = \sum_{j=1}^m w_j\in W$.  We claim that $\bigl|\langle w,
  e_i^*\rangle\bigr|\leqslant 1+\varepsilon$  for each
  $i\in\N$. There are three cases to examine:
  \begin{itemize}
  \item The estimate is obvious for $i<k_1$ because $\langle w,
  e_i^*\rangle =0$ for such~$i$.
  \item Suppose that $i\in[k_j,k_{j+1})$ for some
    $j\in\{1,\ldots,m\}$. Then $\langle w_h,
    e_i^*\rangle =0$ for $h>j$, so
    \[ \bigl|\langle w,
    e_i^*\rangle\bigr|\leqslant \sum_{h=1}^{j-1}\bigl|\langle w_h,
    e_i^*\rangle\bigr| + \bigl|\langle w_j,
    e_i^*\rangle\bigr|\leqslant \frac{(j-1)\varepsilon}{m}+\lVert
    w_j\rVert_X\leqslant \varepsilon+1. \]
  \item  Finally, for $i\geqslant k_{m+1}$, $\bigl|\langle w_j,
    e_i^*\rangle\bigr|\leqslant\varepsilon/m$ for each
    $j\in\{1,\ldots,m\}$, so $\bigl|\langle w,
    e_i^*\rangle\bigr|\leqslant \varepsilon$.
  \end{itemize}
  This establishes the claim, and consequently $\lVert S
  w\rVert_\infty \leqslant1+\varepsilon$
  by~\eqref{formalidXtoc0Eq1}.

  For each $j\in\{1,\ldots,m\}$, we have $\lVert S
  w_j\rVert_\infty\geqslant\varepsilon$, so another application
  of~\eqref{formalidXtoc0Eq1} enables us to choose $h_j\in\N$ such
  that $\bigl|\langle w_j,
  e_{h_j}^*\rangle\bigr|\geqslant\varepsilon$. We note that
  necessarily $h_j\in\left[k_j,k_{j+1}\right)$, and therefore the set
    $\{h_1,h_2,\ldots,h_m\}$ belongs to~$\mathcal{S}_1$ and thus
    to~$\mathcal{S}_n$. This implies that
    \begin{align*}  \lVert w\rVert_X &\geqslant \sum_{i=1}^m
    \bigl|\langle w, e_{h_i}^*\rangle\bigr| =  \sum_{i=1}^m
    \biggl|\Bigl\langle\sum_{j=1}^i w_j,
    e_{h_i}^*\Bigr\rangle\biggr|\\ &\geqslant
    \sum_{i=1}^m
    \biggl(\bigl|\langle w_i, e_{h_i}^*\rangle\bigr| -
    \sum_{j=1}^{i-1}\bigl|\langle w_j,
    e_{h_i}^*\rangle\bigr|\biggr)\geqslant
    \sum_{i=1}^m\biggl(\varepsilon-
    \frac{(i-1)\varepsilon}{m}\biggr) = \frac{(m+1)\varepsilon}{2}.
    \end{align*}
    Combining the above estimates, we conclude that
    \[ 1+\varepsilon\geqslant  \lVert
    Sw\rVert_\infty\geqslant \varepsilon \lVert w\rVert_X\geqslant
    \frac{(m+1)\varepsilon^2}{2}, \] which contradicts that we chose
    $m\geqslant2(1+\varepsilon)/\varepsilon^2$. 
\end{proof}

\begin{proof}[Proof of
    Theorem~{\normalfont{\ref{thmsmallideals}\ref{thmsmallideals2}}}]
As  Odell \cite[p.~694]{odell} observed, the space
  $X=X[\mathcal{S}_n]$ is $c_0$\nobreakdash-sat\-u\-rated (in the
  sense that each of its closed, infinite-dimensional subspaces
  contains an~iso\-mor\-phic copy of~$c_0$) because~$X$ embeds
  into~$C[0,\omega^{\omega^n}]$, which is $c_0$-saturated. Sobczyk's
  Theorem implies that every copy of~$c_0$ in~$X$ is automatically
  complemented, so that~$X$ is sub\-projective, and therefore
  $\mathscr{S}(X) = \mathscr{E}(X)$ by Pfaffenberger's
  result~\cite{Pfaf}.

  Let $S\colon X\to c_0$ be the formal identity operator, as in the
  proof of Lemma~\ref{formalidXtoc0} above. Since~$X$ contains a
  complemented copy of~$c_0$, we can choose operators $U\colon c_0\to
  X$ and $V\colon X\to c_0$ such that $I_{c_0} = VU$.  Then
  $US\in\mathscr{S}(X)\setminus\mathscr{K}(X)$ because
  $S\in\mathscr{S}(X,c_0)$ and $V(US)=S\notin\mathscr{K}(X,c_0)$ by
  Lemma~\ref{formalidXtoc0}, and consequently
  $\mathscr{S}(X)\neq\mathscr{K}(X)$.

Finally, let $Q\in\mathscr{B}(X)$ be a projection whose image is
isomorphic to~$c_0$. Then $Q\notin\mathscr{S}(X)$. However, for each
non-trivial spatial ideal~$\mathscr{I}$, say $\mathscr{I} =
\overline{\langle P_N\rangle}$, where $N\in[\N]$, we can factor~$Q$
through~$X_N$ because~$X$ is $c_0$\nobreakdash-sat\-u\-rated, and
therefore $Q\in\mathscr{I}$. This shows that $Q$ belongs to the
intersection on the left-hand side of~\eqref{thmsmallideals2eq}, and
the conclusion follows.
\end{proof}

\section{Some open questions}\label{sectionOpenQ}
\noindent
Theorem~\ref{thmsmallideals}\ref{thmsmallideals2} and its proof raise some natural
questions. To state them concisely, let $X=X[\mathcal{S}_n]$ for some
$n\in\N$, and denote the closure of the ideal of operators on~$X$
which factor through~$c_0$ by~$\overline{\mathscr{G}}_{c_0}(X)$; in
the notation of the proof of
Theorem~\ref{thmsmallideals}\ref{thmsmallideals2},
$\overline{\mathscr{G}}_{c_0}(X) = \overline{\langle Q\rangle}$, and
the argument given in its last paragraph shows that \[
\overline{\mathscr{G}}_{c_0} (X)\subseteq \bigcap\bigl\{\mathscr{I} :
\mathscr{I}\ \text{is a non-trivial spatial ideal
  of}\ \mathscr{B}(X)\bigr\}. \] However, we do not know whether this
inclusion is proper. We also do not know whether
$\mathscr{S}(X)\subseteq\overline{\mathscr{G}}_{c_0} (X)$.

Another, somewhat less precise, question is as follows. It applies to
both $X=T$ and $X=X[\mathcal{S}_n]$ for
$n\in\N$. Theorem~\ref{manychainsofidealsTsirelson}\ref{manychainsofidealsTsirelson3}
states that~$\mathscr{B}(X)$ contains at least continuum many maximal
ideals, but we do not have an explicit description of a single such
ideal. We know that they cannot be spatial, but is it possible to
describe at least some of these maximal ideals explicitly?

\subsection*{Acknowledgements}
The research presented in this paper was initiated when the third- and
second-named authors visited Washington \& Lee University, VA, in
October 2015 and 2016, respectively, supported by Washington \& Lee
Summer Lenfest grants. It was continued when the first-named author
visited the UK in February/March 2017, supported by a Scheme~2 grant
from the London Mathematical Society. Kania's work has also received
funding from GA\v{C}R project 19-07129Y; RVO 67985840 (Czech
Republic). We gratefully acknowledge this support. Finally, we would like to thank the referee for their careful reading of our paper, especially for suggesting what is now stated as the first part of Lemma~\ref{newlemma210619} as a way of simplifying the proof of the second part of that lemma.

\end{document}